%% file: nonparametric.tex
\documentclass[11pt,fleqn,a4paper]{article}
\pdfoutput=1
\usepackage{amstext}
\usepackage{amsthm}
\usepackage{amssymb,amsmath}
\usepackage{url}
\usepackage{verbatim}
\usepackage{graphicx}
\usepackage[pdfpagemode=UseOutlines ,
plainpages=false,hypertexnames=false ,pdfpagelabels ,hyperindex=true
,colorlinks=true]{hyperref}
	\makeatletter%
	\Hy@breaklinkstrue%
	\makeatother%
\usepackage{color}
\definecolor{darkred}{rgb}{0.5,0,0}
\definecolor{darkgreen}{rgb}{0,0.5,0}
\definecolor{darkblue}{rgb}{0,0,0.5}
\hypersetup{colorlinks ,linkcolor=darkblue ,filecolor=darkgreen
,urlcolor=darkblue ,citecolor=black}
\usepackage{jan-sue-abrev-package}

\RequirePackage{natbib}
\RequirePackage{hypernat}
\renewcommand{\cite}{\citet}
\numberwithin{equation}{section} 

\newtheoremstyle{mysc}
  {3pt}
  {3pt}
  {\it}
  {}
  {\color{darkred}\sc}
  {.}
  {.5em}
  {}

\newtheoremstyle{myex}
  {10pt}
  {10pt}
  {\rm}
  {}
  {\color{darkred}\sc}
  {.}
  {.5em}
  {}

\theoremstyle{mysc}\newtheorem{theorem}{Theorem}[section]
\theoremstyle{mysc}\newtheorem{assumption}{Assumption}[section]
\theoremstyle{mysc}\newtheorem{defin}{Definition}[section]
\theoremstyle{mysc}\newtheorem{lemma}[theorem]{Lemma} 
\theoremstyle{mysc}\newtheorem{corollary}[theorem]{Corollary} 
\theoremstyle{mysc}\newtheorem{proposition}[theorem]{Proposition} 
\theoremstyle{myex}\newtheorem{remark}{Remark}[section]
\theoremstyle{myex}

\makeatletter%
\def\@fnsymbol#1{\ensuremath{\ifcase#1\or 1 \or 2 \or 3\or 4\or  ,\or \star \or
g\or h\or i\else\@ctrerr\fi}}%
\makeatother%

\author{Jan Johannes\thanks{University Heidelberg, 
Institute of Applied Mathematics, Im Neuenheimer Feld, 294, 
D-69120 Heidelberg, Germany, \url{johannes@statlab.uni-heidelberg.de}} 
\and Suhasini Subba Rao\thanks{Texas A\&M University, Department of  
Statistics, College Station, Texas 77843-3143, U.S.A. 
\url{suhasini.subbarao@stat.tamu.edu}}}
\title{Nonparametric estimation for dependent data with an application to 
panel time series \thanks{ 
This work was partially supported by the DFG (DA 187/12-3).}}
\date{\today}

\begin{document}

\maketitle

\begin{abstract}

In this paper we consider nonparametric estimation for dependent data, 
where the observations do not necessarily come from a linear process. 
We study density estimation and also discuss
associated problems in nonparametric regression using the 2-mixing
dependence measure. We compare the results under 2-mixing with 
those derived under the assumption that the process is linear.  

In the context of panel time series where one 
observes data from several individuals, it is often too strong to 
assume the joint linearity of
processes. Instead the methods developed in this paper enable us to 
quantify the dependence through 2-mixing which allows for nonlinearity. 
We propose an
estimator of the panel mean function and obtain its rate of convergence. We
show that under certain conditions the rate of convergence can be improved by
allowing the number of individuals in the panel to increase with time.
\end{abstract}

\begin{tabbing}
\>nonlinear processes, panel time series.\\[.2ex]
\end{tabbing}
\input{introduction}
\input{assumption}
\input{density_estimation}
\input{nonparametric_regression}
\input{panel}
\input{discussion}
\subsection*{Acknowledgements}
The authors are grateful to Professor Rainer Dahlhaus for making 
several useful suggestions. 

\appendix
{
\section{Appendix: Proofs}\label{sec:proofs}
\input{appendix-proofs-nonpar-density}
\input{appendix-proofs-nonpar-regression}
\input{appendix-proofs-nonpar-panel}
\input{appendix-proofs-cov}
\bibliography{nonparametric}
}
\end{document}

%% file: introduction.tex
\section{Introduction}
Nonparametric estimation for dependent observations has a long
history in statistics. \cite{p:ros-70} first studied
density estimation for dependent data. Since then several authors have 
considered
nonparametric estimation under various assumptions. For example, 
\cite{p:hal-har-90}, \cite{p:gir-kou-96}, \cite{p:mie-97} and \cite{p:vie-03}
consider density estimation for linear processes which have long memory, 
whereas \cite{p:che-rob-91} consider density estimation for random variables
which are nonlinear functions of a linear process. 
A notable result, is that they show if the observations 
were from a linear process and have short memory, then the 
usual rate of convergence, known for independent observations, 
also holds for dependent observations. 
On the other hand, for long memory processes, the rate of convergence is different.
Interestingly, despite long memory influencing the rate of convergence, 
there is no influence of long memory on the bandwidth choice, which is 
same regardless of short or long memory. In other
words, if the observations come from a linear process,  
a larger bandwidth does not improve the rate of convergence of the density
estimator. Similar results
can also be derived for nonparametric regression problems (c.f.
\cite{p:hal-har-90b},  \cite{p:che-rob-94} and  \cite{p:cso-95,p:cso-mie-99,p:cso-mie-01}).
However, usually it is assumed that the 
observations come from a linear process or are functions of a linear process.
In the case of linearity, the joint density of the observations can be 
characterised (in some sense) in terms of the autocovariances. It is this
representation that allows for the mean squared error of the nonparametric estimator 
to be derived in terms of the autocovariance function. However this
result does not necessarily hold when the process is nonlinear. 

The assumption of linearity can be relaxed by using the notion of 
2-mixing (see \cite{b:bos-98}), and in this paper we obtain rates of convergence 
for processes which are 2-mixing. 
Unlike the autocovariance function, 2-mixing can be considered as a measure of 
dependence between two random variables (see Definition \ref{defin:size},
below) and the {\it 2-mixing size} quantifies this dependence: 
a large mixing size indicates little dependence, whereas a small mixing size 
indicates large dependence. The 2-mixing size 
can be established for several types of processes, for example, linear
processes, see \cite{p:ath-pan-86},  \cite{p:cli-pu-99}, \cite{p:cha-74} and 
the Appendix \ref{sec:cov-2mix} (noting that 
strong mixing implies 2-mixing, though the converse is not necessarily true) 
and nonlinear processes, see \cite{p:tjo-95}, 
\cite{p:bou-98} and \cite{p:bas-03}. 
Assuming that the 2-mixing size is sufficiently large, 
\cite{b:bos-98} obtains the rate of convergence of several nonparametric estimators.
However despite, there being extensive literature on nonparametric estimation
for linear processes and some on nonparametric estimation for 
processes which are 2-mixing with a sufficiently large 2-mixing size, as far
as we are aware very little exists on nonparametric estimation for nonlinear 
processes whose 2-mixing size is not sufficiently large 
(for example the ARCH$(\infty)$ process, which is a nonlinear process and can
have  a small mixing size). 
In this paper we address this issue, and consider 
nonparametric estimation for dependent data and formulate the results in terms 
of the 2-mixing size. We study both density estimation and also 
nonparametric regression problems. 
A natural application of the methodology proposed in this paper   
is to panel time series, where
one observes several individuals over time and associated with 
each individual are regressors which are known to influence the 
individual. We note that even in the case that an individual 
comes from a linear time series, there is no guarantee that the dependence
between individuals is also linear. Therefore we quantify the dependence in
terms of the 2-mixing size within and between individuals over time, and
consider nonparametric estimation for panel time series within this framework. 
 
In Section \ref{sec:density} we consider kernel density estimation, in particular
we obtain the sampling properties of the 
Rosenblatt-Parzen kernel estimator and obtain a bound for the mean squared error under 
the assumption that the time series are stationary and 2-mixing. We show that, 
like the long memory process, the 2-mixing size can influence the rate of convergence.
But unlike the long memory process, a much larger 2-mixing size is required
to obtain the usual rate of convergence. Moreover, the optimal bandwidths for 
the bounds obtained are influenced by the 2-mixing size - the smaller the 2-mixing 
size the larger the bandwidth.  
We demonstrate that several problems could arise if 
one were to falsely suppose that observations were from a linear process, when
they do not. 
For example, if the usual optimal bandwidth for linear processes 
were used on 
nonlinear processes, the mean squared error may no longer converge to zero. 
Thus our results give a warning to practitioners who
apply well known results for the linear process, 
without checking whether the process  is linear or not.

In Section \ref{sec:non:reg} we consider nonparametric regression for
dependent data. We discuss this with reference to  
two models. First we suppose the response and explanatory
variables $(X_{t},Z_{t})$ satisfy (i) $X_{t} = \varphi(Z_{t}) + 
h(Z_{t})\eta_{t}$, 
where $\{\eta_{t}\}$ and $\{Z_{t}\}$ are independent of each other, and
secondly we assume the conditional expectation satisfies (ii)
$\Ex(X_{t}|Z_{t})  = \varphi(Z_{t})$. We observe that the latter model includes
the former model as a special case. We estimate $\varphi(\cdot)$ using the 
classical kernel estimator and derive rates of convergence 
similar to those obtained for the 
density estimator. But in the
case of model (i) the rate of convergence depends on two factors, the 
2-mixing size of $\{Z_{t}\}$ and 
the rate of decay of the autocovariance function of $\{\eta_{t}\}$, 
whereas for model (ii)
the rate of convergence is determined by the mixing size of the multivariate 
random process $\{(X_{t},Z_{t})\}$. 

Panel time series are often used to model the relationship
and dynamics between several individuals observed over time, and recently
\cite{p:tjo-04} and \cite{p:mam-05} have used nonparametric 
methods in this context.  Typically it is assumed that the dependence between 
individuals is linear, however this assumption is often too strong, as there
could be nonlinear interactions between the individuals. 
In Section \ref{sec:non:pan} we consider estimation for 
nonparametric panel time series, but allow for nonlinear dependence between 
individuals by quantifying their dependence through their 2-mixing sizes. 
Let $X_{t,i}$ denote the observation of the $i$th  
individual at time $t$, where $i=1,\ldots,N$ and $t=1,\ldots,T$. We also 
assume that we observe some explanatory variables $Z_{t,i}$
which influence $X_{t,i}$.
We suppose the influence is common over all
individuals, that is the response and explanatory variables $(X_{t,i},Z_{t,i})$ 
satisfy the relation $\ex[X_{t,i}|Z_{t,i}=z]=\varphi(z)$ for all $i\in\N$ and 
$t\in \Z$. We 
propose a kernel based estimator for $\varphi(\cdot)$ and derive bounds for
the deviation. 
In panel data it is often observed that there is temporal dependence for each 
individual and also dependence
between individuals. We model this by assuming two different 2-mixing sizes.
We show that the rate of convergence of the estimator of $\varphi$ when the 
number of individuals $N$ is kept fixed and $T\rightarrow \infty$, 
is similar, to the rate of convergence of the nonparametric estimator of model (ii)
considered in Section \ref{sec:non:reg}. 
However,  we show that the rates can be improved if we
allow $N$ to increase with $T$. Furthermore,  
if the mixing size is sufficiently large we can obtain the usual nonparametric
rate of convergence obtained for iid random variables.  

All the proofs can be found in the appendix. Also some 2-mixing inequalities for
linear processes used here are included in the appendix.

%% file: assumption.tex
\section{Notation}
In this section we introduce some definitions that will be used in the paper. 
Note we will assume all the necessary densities exist.

We start by defining the multiplicative kernel.
\begin{defin}\label{non:reg:defin:K} 
For all $w = (w_1,\ldots,w_d) \in \R^d$, $K$ is 
a multiplicative kernel (see \cite{b:S92}) of order $r$, i.e. 
$K(w) = \Pi_{i=1}^d \ell(w_i)$ 
where $\ell$ is a univariate, even function such that
        \begin{displaymath}
        \int du \; \ell(u) = 1, \quad \int du \; u^i \ell (u) = 0 
        \end{displaymath}
for all $i=1,\ldots,r-1$ and there exists a  constant 
$S_{K}$ such that 
        \begin{displaymath}
        [\int du \; |u|^r \ell(u)]^d = S_K.
        \end{displaymath}
\end{defin}
Let $K_{b}(z):=b^{-d}K(z/b)$, where $b>0$ is a bandwidth. 
Below we define the smoothness class which we use to 
bound the bias of the estimators. 
\begin{defin}
\label{non:reg:defin:G}
For $s,\triangle> 0$, the space $\mf{G}^{d}_{s,\triangle}$ 
is the class of functions $g:\R^{d}\to \R$ satisfying: $g$ is  everywhere
$(m-1)$-times partially differentiable for $m-1 < s\leqslant m $; where 
for some $\rho  >0$ and for all $x$, the inequality 
        \begin{displaymath}
        \sup_{y:|y-x|<\rho} \frac{\left| g(y)-g(x)-Q(y-x)\right|}{|y -x|^{s}}
\leqslant \triangle, 
\end{displaymath} 
holds true with $Q=0$ when $m=1$ and for $m>1$, $Q$ is an 
$(m-1)$th-degree homogeneous polynomial in $y-x$, whose coefficients are
the partial derivatives of $g$ of orders $1$ to $m-1$ evaluated at 
$x$;  and $\triangle$ is a finite constant.
\end{defin}
For brevity, we use the standard notation $\wedge$ to denote minimum and 
$\vee$ to denote maximum.

%% file: density_estimation.tex
\section{Kernel density estimation}\label{sec:density}
Suppose we observe the stationary time series
$\{Z_{1},\ldots,Z_{T}\}$, and let $f$ denote the density of $Z_{t}$.
The most popular estimator of $f$, is the Rosenblatt-Parzen kernel estimator 
\begin{eqnarray}
\label{eq:density_estimator}
\hat{f}(u) = \frac{1}{T}\sum_{t=1}^{T}K_{b}(Z_{t}-u),  
\end{eqnarray}
where $K_{b}(z) = b^{-1}K(\frac{z}{b})$, $b>0$ is a bandwidth, and $K$
is a multiplicative kernel (see Definition \ref{non:reg:defin:K}). 
In this section we investigate the sampling properties of the kernel 
density estimator defined above. The 
dependence of the process $\{Z_{t}\}$ is quantified 
in terms of its 2-mixing size.
\begin{defin}\label{defin:size}
\begin{itemize}
\item[(i)] A process $\{Y_{t}\}$ is said to be 2-mixing with size
$\mathfrak{v}$ if for all $t\ne\tau$ 
\begin{displaymath}
\sup_{A\in \sigma(Y_{t}),B\in
\sigma(Y_{\tau})} |P(A\cap B) - P(A)P(B)|\leq C|t-\tau|^{-\mathfrak{v}}.
\end{displaymath}
for some $C<\infty$ independent of $t$ and $\tau$.
\item[(ii)] The covariance of a stationary process $\{Y_{t}\}$ has
size $\mathfrak{u}$ if for all $t\ne\tau$, 
$|\cov(Y_{t},Y_{\tau})|\leq C|t-\tau|^{-\mathfrak{u}}$ for some $C<\infty$ 
independent of $t$ and $\tau$.
\end{itemize}
\end{defin}
We note that the covariance is a measure of linear dependence, whereas 
2-mixing is a generalisation of this, and can be considered as a 
measure of dependence. 2-mixing is quite a general notion, which is satisfied
by several processes. For example, under certain conditions on the innovations, 
most linear models are 2-mixing (see Appendix \ref{sec:cov-2mix}, and 
\cite{p:ath-pan-86},  
\cite{p:cli-pu-99} and \cite{p:cha-74}, where strong mixing is shown). 
Further, under additional conditions
on the innovations and the parameters, ARCH/GARCH processes are also strongly mixing 
(c.f. \cite{p:tjo-95}, 
\cite{p:bou-98} and  \cite{p:bas-03}) which implies that they also
2-mixing. Most of the 
results and bounds in this paper are derived using 2-mixing. In general, 
the larger 
the mixing size the faster the rate of convergence. For example,  in the case of iid observations 
(the 2-mixing size can be treated as $\infty$) 
using just a few observations, information over the entire domain 
of the density 
function can be obtained. On the other hand,  
a sample which has a small mixing size (so tends to be clustered about certain points) 
will require a much larger number of observations to give the same information.

We first derive a bound for the  mean 
squared error (MSE) $\Ex|\hat f(z) - f(z)|^{2}$ using 
only minimal assumptions on the distribution of $\{Z_{t}\}$.
\begin{proposition}\label{prop:cov:bra}
Suppose the stationary process $\{Z_{t}\}$ is 2-mixing with size 
$\gv$ and the marginal density $f$ of $Z_{t}$  and its second derivative $f^{\prime\prime}$ are both
uniformly bounded.  Let $\hat{f}$ be defined as in
(\ref{eq:density_estimator}), where $K$ is a rectangular kernel, i.e.,  
$K(x) = 1$ if $x\in [-1/2,1/2]$ and zero otherwise. Then we have
\begin{eqnarray*}
\Ex|\hat f(z) - f(z)|^{2} = 
O(b^{4} + T^{-[\gv\wedge 1]}b^{\frac{-[(\gv \vee 1) + 1]}{\gv \vee 1}})  = 
\left\{
\begin{array}{cc}
O(b^{4} + T^{-1}b^{-\frac{\gv+1}{\gv}}) & \gv > 1 \\
O(b^{4} + T^{-\gv}b^{-2}) & \gv \leq 1\\
\end{array}
\right.
\end{eqnarray*}
\end{proposition}
\noindent\textcolor{darkred}{\sc Proof.} To prove the result we will bound the risk using the standard variance bias decomposition. First the bias: as we are using a rectangular kernel and 
$f^{\prime\prime}$ is uniformly bounded, it is clear that 
$\Ex(\hat f(z)) = f(z) + O(b^{2})$.
To obtain a bound for the variance we require a bound for the 
covariances inside the variance expansion 
$T^{2}\cdot\var(\hat f(z)) =
\sum_{t,\tau}\cov[K_{b}(Z_{t}-z),K_{b}(Z_{\tau}-z)]$.
Since $\{Z_{t}\}$ is 2-mixing with size 
$\gv$ by using the covariance inequality in \cite{p:bra-96} 
(see also \cite{p:rio-93}) we have
\begin{multline}\label{eq:cov:bra}
|\cov[K_{b}(Z_{t}-z),K_{b}(Z_{\tau}-z)]| \\
\hfill\leq 4\cdot \int_{0}^{\infty}
\int_{0}^{\infty}\min\left(C|t-\tau|^{-\gv},
P(|K_{b}(Z_{t}-z)|>x), P(|K_{b}(Z_{\tau}-z)|>y)\right)dxdy.
\end{multline}
Studying $P(|K_{b}(Z_{t}-z)|>x)$ and recalling that $K(\cdot)$ is a rectangular
kernel we can show that
\begin{align*}
 P(|K_{b}(Z_{t}-z)|>x) &= 
\left\{
\begin{array}{cc}
 0, & \textrm{ if } x> 1/b; \\
 P(Z_{t} \in [z-b/2,z+b/2]), & \textrm{ otherwise.}
\end{array}
\right.
\end{align*}
By using the mean value theorem we have 
$P(X_{t} \in [z-b/2,z+b/2]) = bf(\tilde{z})$, for some $\tilde{z}\in [z-b/2,z+b/2]$.
Substituting this into (\ref{eq:cov:bra}) leads to 
\begin{align}\nonumber
|\cov[K_{b}(Z_{t}-z),K_{b}(Z_{\tau}-z)]| &\leq 4\cdot \int_{0}^{1/b}
\int_{0}^{1/b}\min\left(C|t-\tau|^{-\gv},
 b\cdot f(\tilde{z})\right)dxdy \\\label{eq:cov:bra:e1}
&= 4\cdot  b^{-2}\min\left(C|t-\tau|^{-\gv},
 b\cdot f(\tilde{z})\right).
\end{align}
Altogether this yields the bound 
\begin{eqnarray*}
T^{2}\cdot\var(\hat f(z)) \leq 4\sum_{t,\tau}b^{-2}\min\left(C|t-\tau|^{-\gv},
 b\cdot f(\tilde{z})\right).
\end{eqnarray*} 
Examining the minimum inside the summand above, we partition the sum into two
parts which we bound separately (for the details see the proof of 
Theorem \ref{non:theo}, in the Appendix). 
Finally recalling that 
$[\Ex(\hat f(z)) - f(z)]^{2} = O(b^{4})$ leads to the desired result. \hfill $\Box$

\vspace{3mm}
We observe, in the proof above, that besides the 2-mixing condition we 
do not have any assumptions on the joint distribution of $(Z_{t},Z_{\tau})$.
The cost of using such weak assumptions is that the usual bound $O(b^{4} + (bT)^{-1})$ for the MSE, obtained for independent observations, is not 
achieved. Even for large $\gv$ the 2-mixing size has an influence on the
bound. However,  introducing some assumptions on the joint densities of $\{Z_{t}\}$ 
allows us to tighten the bound derived in (\ref{eq:cov:bra:e1}) and, hence for
a sufficiently large mixing size 
$\gv$ to recover the usual  bound $O(b^{4} + (bT)^{-1})$ for the MSE (we note
that the rest of the proofs
in this section and the subsequent sections require more subtle arguments,
and these can be found in the appendix). 
\begin{assumption}[Densities and kernels]\label{dens:est:ass:tec-part}
\begin{itemize}
\item[(i)] The marginal density $f$ is uniformly bounded.
\item[(ii)] For each $t,\tau\in\Z$ let 
$f^{(t,\tau)}$ denote the joint density of $(Z_{t},Z_{\tau})$. 
 Define\footnote{We use the notation $f\otimes g (x,y) = f(x)g(y)$ and 
$\|f\|_{p}=(\int |f(x)|^{p}dx)^{1/p}$.}
$F^{(t,\tau)} :=  f^{(t,\tau)} - f\otimes f$. 
Then $\|F^{(0,t)}\|_{p_{F}}$ is uniformly bounded in $t$ for some
$p_{F}>2$  and 
we define $q_{F}=1-2/p_{F}$.
\item[(iii)] The univariate kernel $K$ is uniformly bounded and has a finite 
first and second moment, i.e.,
$\|K\|_{1}<\infty$ and $\|K\|_{2}<\infty$. 
\end{itemize}
\end{assumption}

We use these assumptions to derive a bound for the MSE of  
the density estimator. 
\begin{theorem}\label{non:theo}
Let us suppose the stationary time series $\{Z_{t}\}$ is $2$-mixing with size 
$\gv$ and 
Assumption \ref{dens:est:ass:tec-part} is fulfilled for some $q_{F}\in(0,1)$. 

Let $\hat f$ be defined as in (\ref{eq:density_estimator}), where
$K$ is a univariate kernel of order $r$.
In addition assume that $f\in \mf{G}^{1}_{s,\triangle}$ 
for some $\triangle,s>0$ (see Definition \ref{non:reg:defin:G}), 
and let $\rho=r\wedge s$. Then we have for all $z\in\R$
\begin{equation*}
\Ex|\hat f(z)- f(z)|^{2}=O\Bigl(b^{2\rho} + b^{-1}\cdot T^{-1} + 
b^{-2-q_{F}(1-[\gv\vee 1])}\cdot T^{-[\gv\wedge 1]}  \Bigr), \quad T\to\infty.
\end{equation*}
\end{theorem}
For ease of presentation we have only stated the result for univariate $\{Z_{t}\}$,
however it is straightforward to extend this result for multivariate $\{Z_{t}\}$. 
Indeed, the proof of the theorem given in the Appendix is derived for
random vector $\{Z_{t}\}$ (as we require the multivariate case 
in Section \ref{sec:non:reg}).

\begin{remark}\label{remark:density}
We note that in the bound given in Theorem \ref{non:theo} the second term  dominates  
the third term when $\mathfrak{v} > 1+ 1/q_{F}$. Conversely, when 
$\mathfrak{v} < 1+ 1/q_{F}$
the third term dominates the second term. Moreover, the third term can be partitioned
into two further cases, when $1<\mathfrak{v}\leq 1+ 1/q_{F}$ and when 
$\mathfrak{v}\leq 1$. 
This means that  Theorem \ref{non:theo} can be written as 
\begin{itemize}
\item[(i)] if $\mathfrak{v} > 1+ 1/q_{F}$, then  
$\Ex|\hat f(z)- f(z)|^{2}=O\Bigl(b^{2\rho} + \frac{1}{bT} \Bigr)$;
\item[(ii)] if  $1<\mathfrak{v}\leq 1+ 1/q_{F}$, then 
$\Ex|\hat f(z)- f(z)|^{2}=O\Bigl(b^{2\rho} + 
\frac{1}{b^{(2+q_{F}(1-\mathfrak{v}))}T} \Bigr)$;
\item[(iii)] if  $\mathfrak{v}\leq 1$ then 
$\Ex|\hat f(z)- f(z)|^{2}=O\Bigl(b^{2\rho} +  \frac{1}{b^{2}T^{\mathfrak{v}}}\Bigr)$;
\end{itemize}
as $T\rightarrow \infty$. 

Studying the three bounds, we see that the bound increases linearly with 
$\gv$ for $0 \leq  \mathfrak{v} \leq 1$, 
after this point there is a change in behaviour and the increase is more
gradual. The bound plateaux when $\mathfrak{v} > 1+1/q_{F}$, after this point 
we have the usual nonparametric bound obtained for iid observations.
There is also a continuity in the three bounds. More precisely, when
 $\mathfrak{v}$ is at the boundary of $1$ and $1+q_{F}^{-1}$, there is a 
continuous transition between the bounds. 
\hfill $\Box$
\end{remark}

We now consider the rate of convergence using the optimal bandwidth $b^{*}$.
\begin{corollary}\label{cor:theo}
Suppose the assumptions in Theorem \ref{non:theo} are satisfied and $r\geq s$. Let 
$b^*\approx T^{-\gamma/(2s+1)}$ with 
\begin{equation}\label{cor:theo:gamma}
\gamma:=\left\{\begin{array}{ll}
1,&\gv>1+1/q_{F};\\\,
[\gv\wedge 1] \cdot \frac{2s +1}{2s+(2+q_{F}(1-[\gv\vee 1]))},&1+1/q_{F}\geq \gv.
\end{array}\right.
 \end{equation}
Then for all $z\in\R$ we have  $ \Ex|\hat f(z)- f(z)|^{2} = 
O\Bigl(
T^{-\frac{2s}{2s+ 1}\cdot \gamma}\Bigr)$ as $T\rightarrow \infty$.
\end{corollary}
In other words, if  $b^*\approx T^{-\gamma / (2s+1)}$, then we have  
\begin{equation}\label{non:theo:cor}
\Ex|\hat f(z)- f(z)|^{2}:=\left\{\begin{array}{ll}
O\Bigl(T^{-\frac{2s}{2s+ 1}}\Bigr) ,& \gv>1+1/q_{F};\\
O\Bigl(T^{-\frac{2s}{2s+1}\cdot(\frac{2s +1}{2s+(2+q_{F}(1-\gv))})}\Bigr),
&1+1/q_{F}\geq \gv>1;\\
O\Bigl(T^{(\gv\cdot\frac{2s +1}{2s+2})}\Bigr),&1\geq \gv.
\end{array}\right..
 \end{equation}

We note that if 
$\sup_{z}|f(z)|<\infty$ and $\sup_{t}\sup_{z}|f^{(0,t)}(z)|<\infty$ (both the 
density and the joint densities are uniformly bounded),  then 
uniformly in all $t$, $\|F^{(0,t)}\|_{\infty}<\infty$. This means $q_{F} = 1$, and the
bound can be divided into the three cases where 
$\mathfrak{v}\leq 1$, $ 1\leq \mathfrak{v}\leq 2$  and $\mathfrak{v}\geq 2$.
On the other hand when $\|F^{(0,t)}\|_{p_{F}}<\infty$ for only a finite
$p_{F}$, then 
$q_{F}<1$ and $\mathfrak{v} >1+q_{F}^{-1}> 2$ to be sure of the usual nonparametric
bound.

Referring to Corollary \ref{cor:theo}, we observe that when 
$\mathfrak{v}<1+q_{F}^{-1}$, then the
optimal bandwidth $b^{*}$ is much larger than usual optimal bandwidths
encountered in nonparametric regression 
($b\approx T^{\frac{-1}{2s+1}})$. We discuss this further in Section
\ref{sec:density-nonlinear}.

\subsection{A comparison of the MSEs for linear processes}\label{sec:linear}
In this section we compare the MSE in Theorem \ref{non:theo} with the 
results obtained under the stronger condition that the observations $\{Z_{t}\}$ come 
from a linear process. We will use the results in Appendix \ref{sec:cov-2mix} and 
show that if the process were linear, and not just 
mixing, that then the rate of convergence is better than the rate obtained in Corollary
\ref{cor:theo}. However, in Section \ref{sec:density-nonlinear} 
we demonstrate that by misspecifying the process to be 
linear, can lead to several problems with the density estimator, including 
bounds which do not converge to zero.  

\begin{figure}[h]
\begin{center}
\includegraphics[height=10cm,width=10cm]{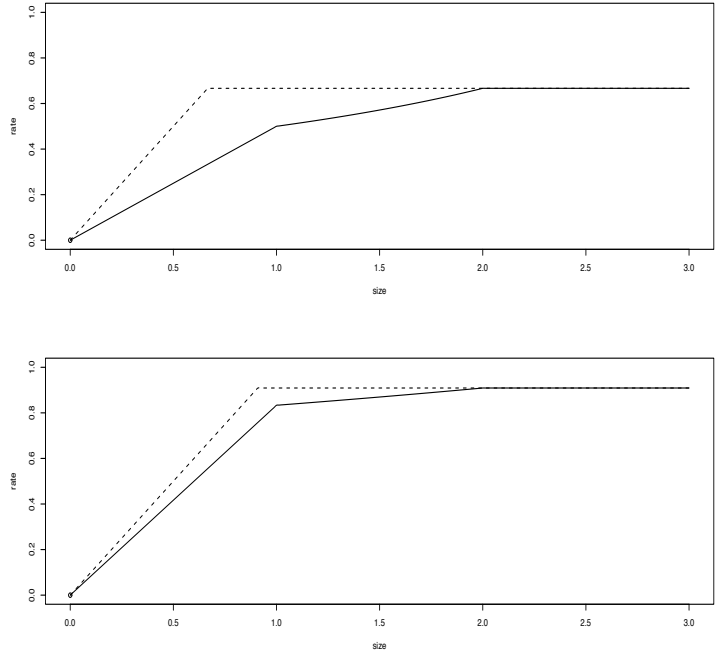}
\end{center}
\caption{The top and bottom plot corresponds to 
$\rho =  (r\wedge s) = 1$ and  
$\rho =(r\wedge s)= 5$, respectively.  
The $x$-axis is the covariance and mixing size (assuming
both are the same) and the $y$-axis is the indice $\delta$ in the MSE
$\Ex|\hat{f}(x)-f(x)|^{2} = O(T^{-\delta})$. 
The solid line is the MSE using 2-mixing and dotted line is the 
MSE when $\{Z_{t}\}_{t}$ is a linear process. 
We have assumed that $q_{F} = 1$ (in other words $\|F^{(0,t)}\|_{\infty}<\infty$).
\label{fig:1}}
\end{figure}

Let us suppose $\{Z_{t}\}$ has a linear process representation and 
satisfies 
\begin{eqnarray}\label{eq:MAinfty}
Z_{t} = \sum_{j=0}^{\infty}a_{j}\varepsilon_{t-j},
\end{eqnarray}
where the innovations $\{\varepsilon_{t}\}$ are iid.
Under the assumptions in Lemma \ref{lemma:MAdensity} (see the Appendix), it
can be shown that 
$\cov(K_{b}(Z_{0}),K_{b}(Z_{t}))=O(\cov(Z_{0},Z_{t}))$. Using this as the basis,  
\cite{p:hal-har-90}, \cite{p:gir-kou-96}, \cite{p:mie-97} and \cite{p:vie-03}) 
have shown that 
if the  kernel is of order $r\geq s$ (where $r$ is the order of
the kernel and $s$ is the smoothness of $f$, see Definitions 
\ref{non:reg:defin:K} and \ref{non:reg:defin:G}),  the MSE is 
\begin{eqnarray}\label{eq:linear-rates1}
\Ex|\hat f(z)- f(z)|^{2} = O\left(b^{2s} + 
\frac{1}{bT} + \frac{1}{T}R_{T}\right),
\end{eqnarray}
where $R_{T} = \sum_{t=1}^{T}|\cov(Z_{0},Z_{t})|$. 
It is clear that both 
$\cov(Z_{0},Z_{t})$ and $R_{T}$ depend on the rate of decay of the 
parameters $a_{j}$. We observe if
$|a_{j}|\leq Cj^{-\theta}$, then 
\begin{eqnarray*}
\begin{array}{lllll}
|\cov(Z_{0},Z_{t})|=O(T^{-2\theta + 1}) &\textrm{ and }&
R_{T} =  O((\log T)T^{-(2\theta - 1)+1}) & \textrm{if} &
1/2 < \theta \leq  1 \\
|\cov(Z_{0},Z_{t})| = O(T^{-\theta}) &\textrm{ and }&
R_{T}=O(T^{-\theta +1 }) & \textrm{if }& \theta > 1. 
\end{array}
\end{eqnarray*}
Substituting these rates into (\ref{eq:linear-rates1}) we see that the bound of the MSE depends on $\theta$. 
We recall that a process $\{Z_{t}\}$ is called a short memory process if 
$\sum_{t}|\cov(Z_{0},Z_{t})| < \infty$, otherwise it is called a long memory process.
Now studying (\ref{eq:linear-rates1}) we see that  
$R_{T}$ does not depend on the bandwidth $b$. In other words 
 long memory has {\it no influence} on the choice of the optimal bandwidth. 
To summarise, the rate of convergence for observations coming from a linear process is 
\begin{eqnarray}\label{eq:linear-rates}
\Ex|\hat f(z)- f(z)|^{2} \leq \left\{
\begin{array}{lll}
O(T^{-(2\theta -1 )}), & \textrm{ if }& 2\theta -  1 \leq  \frac{2s}{2s+1}; \\
O(T^{\frac{-2s}{2s+1}}), & \textrm{ if }& 2\theta - 1  >  \frac{2s}{2s+1}.
\end{array}
\right.
\end{eqnarray}
We now compare these results to (\ref{non:theo:cor}), in particular when 
$\mathfrak{v} \leq 1$, we have 
\begin{eqnarray}\label{eq:linear-bad-rates}
\Ex|\hat f(z)- f(z)|^{2} = O(T^{- \mathfrak{v}  + \frac{2\mathfrak{v}}{2s+2}}). 
\end{eqnarray}
It is difficult to directly compare (\ref{eq:linear-rates}) and 
(\ref{eq:linear-bad-rates}), 
since (\ref{eq:linear-rates}) is in terms of its long memory parameter whereas
(\ref{eq:linear-bad-rates}) is in terms of its mixing size $\mathfrak{v}$. 
However in the
special case that $\{Z_{t}\}$ is Gaussian (and thus linear), 
there is a one-to-one correspondence, for example, 
if $ 2\theta -  1 \leq 1 $ then  the covariance size and mixing size are the same, and 
$\mathfrak{v} = (2\theta -  1)$. Noting that the Gaussian density is  analytic, 
the rate of convergence is determined by the order of the kernel $r$. In this case,  
the rates in (\ref{eq:linear-rates}) are better than those in 
(\ref{eq:linear-bad-rates}), 
but as the order $r$ increases the two rates become close. 
We illustrate the case when the mixing and the covariance sizes are the same
in Figure \ref{fig:1} (for both large and small $r\wedge s$). 
In the non-Gaussian case, where the 2-mixing and covariance size do not necessarily
coincide, $\mathfrak{v} \neq (2\theta -  1)$, we have that 
$(2\theta -  1)\frac{\ell}{2(\ell+1)} \leq \mathfrak{v} \leq  (2\theta -  1) 
\frac{\ell}{(\ell-2)}$, where  
the innovations satisfy $\Ex(|\varepsilon_{0}|^{\ell})  < \infty$
(see (\ref{eq:up:lo:bd}) in the appendix). 
In this case it is not clear which rate (\ref{eq:linear-rates}) or 
(\ref{eq:linear-bad-rates}) is better. 
However substituting the lower bound 
$\mathfrak{v}\geq (2\theta -  1)\frac{\ell}{2(\ell+1)}$ into 
Corollary \ref{cor:theo} yields a rate which is less than (\ref{eq:linear-rates}).
In summary better rates of convergence can often be obtained if the observations come from a linear process.
On the other hand,  2-mixing  is a {\it weaker condition, that 
is satisfied by a far wider class of processes}. We consider below the MSE for 
processes which are not linear, and show that misspecifying the model, and assuming
linearity, when the process is nonlinear could severely affect the MSE. 

\subsection{The MSE for nonlinear processes}\label{sec:density-nonlinear}

As far as we are aware, theory is required to bridge the gap for 
processes which are nonlinear but have a small mixing size.
One of the main aims of Theorem \ref{non:theo} is to fill in the gap in the 
theory, and to derive a bound for the MSE when the observations come from  
nonlinear processes with small 2-mixing size.

The joint densities of processes which are nonlinear do not necessarily 
satisfy 
the density decomposition in Lemma \ref{lemma:MAdensity}. 
Without this result it cannot be shown that\linebreak $\cov(K_{b}(Z_{0}),K_{b}(Z_{t}))=
O(\cov(Z_{0},Z_{t}))$, and the rates in (\ref{eq:linear-rates}) do not
necessarily hold. 
Instead, to prove the results, under Assumption \ref{dens:est:ass:tec-part},
we use classical mixing inequalities to tighten the bound given in 
(\ref{eq:cov:bra:e1}) (see the proof of Proposition \ref{prop:cov:bra}).
More precisely, to prove Theorem \ref{non:theo} we show that 
\begin{eqnarray*}
|\cov(K_{b}(Z_{t}-z),K_{b}(Z_{\tau}-z))|\leq C\cdot b^{-2}
\min\left(|t-\tau|^{-\mathfrak{v}},
  b^{(1+q_{F})}\right),
\end{eqnarray*}
where $C$ is a finite constant (see Lemma \ref{a:p:non:dens:l:1}, for the
proof).  

Looking at some of the implications of 
Theorem \ref{non:theo}, we demonstrate below that several problems could arise if 
one were to falsely suppose that the observations come from a linear process, when
they do not.
\begin{itemize}
\item[(i)] In the case of linear processes, the optimal bandwidth has the
same order as the 
optimal bandwidth for iid random variables (regardless of 
long memory). The same is not necessarily true when all that 
is known is that the process is 2-mixing. Moreover, if the mixing size  satisfies
$\mathfrak{v} \leq 1$ and the bandwidth is such that 
$b^{2}T^{\mathfrak{v}}<\infty$, then we  see from Theorem \ref{non:theo}
that the bound does not converge to zero. An important example, is when the `usual' optimal bandwidth for linear or iid data is used (that is $b \approx T^{-\frac{1}{2s+1}}$). In this case, substituting 
$b \approx T^{-\frac{1}{2s+1}}$ into  Theorem \ref{non:theo} leads to the result
\begin{eqnarray*}
\Ex|\hat f(z)- f(z)|^{2} = 
\left\{
\begin{array}{cc}
O(T^{\frac{-2s}{2s+1}}) &  \mathfrak{v}  > 1+1/q_{F} \\
O( T^{\frac{1+q_{F}(1-\mathfrak{v})-2s}{2s+1}})
& 1 <  \mathfrak{v}  \leq 1+1/q_{F} \\
 O(T^{\frac{2 - \mathfrak{v}(2s+1)}{2s+1}}) 
& 0 \leq   \mathfrak{v}  \leq 1 \\
\end{array}
\right.
\end{eqnarray*}
Studying the rates above we see when $\mathfrak{v}  < 1+1/q_{F}$, the rates
are lower than the rates given using the optimal bandwidth (compare 
the above with the rates in Corollary \ref{cor:theo}). 
Moreover, in the case that $\mathfrak{v}\leq \frac{2}{2s+1}$, the 
bound cannot be used to show consistency of the estimator - since the bound
does not even converge to zero. 

In short, to estimate the density at any given
point,  the number of observations (approximately $bT$)
needs to be much larger than in the iid case. 
\item[(ii)] Rather surprisingly even when 
$\sum_{j=1}^{\infty}|\cov(Z_{0},Z_{j})|<\infty$,
the `usual $O(T^{-\frac{2s}{2s+1}})$' rate, may not hold, unlike for linear processes.
However, the usual rate does hold when 
$\mathfrak{v}\geq 1+1/q_{F} > 2$. Therefore, even 
when the mixing and covariance size
are the same, a far larger mixing size may be require to obtain the 
`usual $O(T^{-\frac{2s}{2s+1}})$' rate of convergence.
\end{itemize}

Our results give a cautionary warning to practitioners who
apply the optimal bandwidths for linear processes to 
nonlinear process. In the subsequent sections, where
we consider nonparametric regression problems, the assumptions and proofs
will be more involved, however the underlying message is the same. 
That is, more than just the second order autocovariance function may
have influence on the rate of convergence, 
and the rate of convergence can be severely compromised
if the usual bandwidths were used.

\begin{remark}[Example]
It is almost impossible to estimate the 2-mixing size from the observations, 
in contrast to long memory 
(c.f. \cite{p:gew-por-83}, \cite{p:kue-87} and \cite{p:rob-95}).
However to conclude this section we give an example of a nonlinear process
whose 2-mixing size is less than $1+\delta$, for some $\delta>0$.
Let us consider the ARCH$(\infty)$ process (see \cite{p:rob-91}), where 
$\{Z_{t}\}$ satisfies
\begin{eqnarray*}
Z_{t} = \sigma_{t}\varepsilon_{t} \qquad \sigma_{t}^{2} = a_{0} +
\sum_{j=1}^{\infty}a_{j}Z_{t-j}^{2},  
\end{eqnarray*}
with $\Ex(\varepsilon_{t})=0$ (estimation of ARCH$(\infty)$ parameters
is considered in \cite{p:sub-07}). \cite{p:gir-00} have shown that if for large $t$, 
$a_{t}\approx  t^{-(1+\delta)}$ (for some $\delta<0$) and 
$[\Ex(\varepsilon_{t}^{4})]^{1/2}\sum_{j=1}^{\infty}a_{j} <1$, 
 then $|\cov(Z_{0}^{2},Z_{t}^{2})|\approx  t^{-(1+\delta)}$.
That is, the absolute sum of the covariances is finite, but `only just' if
$\delta$ is small. 
Furthermore, if we assume that $|\varepsilon_{t}| < 1$, 
then it is straightforward 
to show that $Z_{t}$ is a bounded random variable. This means that 
using the mixing inequality for bounded random variables 
(see \cite{b:hal-hey-80}) we can show
\begin{eqnarray*}
 |\cov(Z_{0}^{2},Z_{t}^{2})| \leq   C
\sup_{A\in \sigma(Z_{0}),B\in\sigma(Z_{t})} |P(A\cap B) - P(A)P(B)|,
\end{eqnarray*}
for some $C<\infty$. Altogether this implies that an upper bound for  
the 2-mixing size of the ARCH$(\infty)$ process with 
$a_{t}\approx t^{-(1+\delta)}$, is $\gv \leq (1+\delta)$.  

In other words the 2-mixing size for some ARCH$(\infty)$ process is small, and 
far from the geometric rate often assumed in nonparametric estimation.
A lower bound for the 2-mixing size can be found in \cite{p:sub-07b}. 
\hfill $\Box$\end{remark}

%% file: nonparametric_regression.tex
\section{Nonparametric regression}\label{sec:non:reg}
In this section we consider nonparametric regression, 
with random design, where the observations are dependent. 
It is worth mentioning that there has been extensive research done on  
nonparametric regression with fixed design and dependent errors
(c.f. \cite{p:hal-har-90b}, \cite{p:cso-95}, 
and the references therein). In this case typically,  
one observes $Y_{t}$, where $Y_{t} = \varphi(\frac{t}{T}) + \varepsilon_{t}$
and $\{\varepsilon_{t}\}_{t}$ are stationary random variables with varying degrees of 
dependence. It has been shown that the rate of convergence depends on the covariance 
of $\{\varepsilon_{t}\}_{t}$, in particular their absolute sum,
$\sum_{t=1}^{\infty}|\cov(\varepsilon_{0},\varepsilon_{t})|$. 

In the random design model, one observes the stationary $(1+d)$-dimensional 
 vector time series $\{(X_{t},Z_{t})\}_{t}$, where 
\begin{eqnarray}\label{non:reg:model}
X_{t} = \varphi(Z_{t}) + \varepsilon_{t}
\end{eqnarray}
with $\Ex(X_{t}|Z_{t}=z)=\varphi(z)$ and $\varepsilon_{t} = X_{t} - \Ex(X_{t}|Z_{t})$.
The randomness in this model is determined by two factors: the design 
$\{Z_{t}\}$ and the errors $\{\varepsilon_{t} = X_{t} - \Ex(X_{t}|Z_{t})\}$. 
Therefore, unlike the fixed design model, the rate of convergence of any estimator
of $ \varphi$ must depend on the sampling properties of the 
design density estimator. Thus, it is clear that similar results to those in 
Section \ref{sec:density} should also apply to an estimator of $\varphi$. 

We now define the classical Nadaraya-Watson estimator of 
$\varphi(\cdot)$ and study its sampling properties, under various assumptions on 
$\{(X_{t},Z_{t})\}$.  Let $p(x,z)$ be the joint density of $(X_{t},Z_{t})$. The 
estimator is 
\begin{equation}\label{non:reg:est:classical}
\hat\varphi(z)=\frac{
\hat g(z)}{\hat f(z)},
\end{equation}
where $\hat g(z):=\frac{1}{T}\sum_{t=1}^T X_{t}K_{b}(Z_{t}-z)$
and $\hat f(z):=\frac{1}{T}\sum_{t=1}^TK_{b}(Z_{t}-z)$
are estimators of $g(z) = \int x p(x,z)dx$ and $f(z)$, which is the density of 
$Z_{t}$.

We first consider the sampling properties for a particular class of models
which satisfy (\ref{non:reg:model}). 
Suppose the vector time series $\{(X_{t},Z_{t})\}$ satisfies the representation
\begin{eqnarray}\label{non:reg:model1}
X_{t}= \varphi(Z_{t}) + h(Z_{t})\eta_{t}
\end{eqnarray}
for some $h:\R^d\to \R^+$,  where the time series  $\{Z_{t}\}$ and $\{\eta_{t}\}$ are 
independent of each other. This class of models is similar to 
the fixed design model $X_{t} = \varphi(\frac{t}{T}) + \eta_{t}$, but in (\ref{non:reg:model1}) the 
design is random and
the conditional variance $\var(X_{t}|Z_{t})=h(Z_{t})^{2}\var(\eta_{t})$,
depends on the design. 
This model arises in various applications and we consider one application in Remark 
\ref{remark:example}. We will show in the theorem below that the rate of
convergence depends 
both on the mixing size of the design  $\{Z_{t}\}$, but also on the size of
the covariances of the process $\{\eta_{t}\}$ (which we denote by $\gu$, see 
Definition \ref{defin:size}). 

We require the following assumptions. 
\begin{assumption}[Densities, moments and kernels]\label{non:reg:ass:tec-part1}
\begin{itemize}
\item[(i)]  For some $p>2$ the functions $h^2\cdot f$ and  $|\varphi|^p\cdot
  f$  are uniformly  bounded
and we define  $q:=1-2/p$. 
\item[(ii)]Let 
$f^{(t,\tau)}$ and $F^{(t,\tau)}$ be defined as in 
Assumption \ref{dens:est:ass:tec-part} (ii),
\begin{displaymath}
g^{(t,\tau)}(z_{1},z_{2}):= \ex[X_{t}X_{\tau}| Z_{t}=z_{1},
Z_{\tau}=z_{2} ]\cdot f^{(t,\tau)}(z_{1},z_{2}).
\end{displaymath}
and 
$G^{(t,\tau)} := g^{(t,\tau)} - g \otimes g$. 
Then $\|F^{(t,\tau)}\|_{p_{F}}$ and $\|G^{(t,\tau)}\|_{p_{G}}$ are uniformly 
bounded in 
$t$ and $\tau$ for some $p_{F},p_{G}>2$. We define $q_{F}:=1-2/p_{F}$, 
$q_{G}:=1-2/p_{G}$ 
and $q_{FG}:=q_{F}\wedge q_{G}$.
\item[(iii)] The multiplicative kernel $K$ has finite first and $p$-th moment.
\end{itemize}
\end{assumption}

Studying Assumption \ref{non:reg:ass:tec-part1}(i), we see that it allows for 
various types of growth of the regression function $\varphi$ and the 
conditional variance $h$.
The type of growth depends on the rate the density $f$ decays to zero. 
For example, if $f$ were the Gaussian density, then exponential growth of $\varphi$ 
and $h$ is possible.  However, as we shall demonstrate in the theorem below,
the larger the $p$, such that $\sup_{x}h(x)^p\cdot f(x)<\infty$ and 
$\sup_{x}|\varphi(x)|^p\cdot f(x)<\infty$, then the faster the rate of convergence of 
$|\hat \varphi(z)- \varphi(z)|^{2}$.

\begin{theorem}\label{non:reg:theo-model1}
Suppose the stationary time series 
$\{(X_{t},Z_{t})\}$ satisfies (\ref{non:reg:model1}),  
$\{Z_{t}\}$ is $2$-mixing with size $\gv$  and the autocovariance of the time series 
$\{\eta_{t}\}$ has size $\gu$.  
Suppose Assumption \ref{non:reg:ass:tec-part1} is fulfilled for some 
$q,q_{FG}\in(0,1)$. 

Let the estimator $\hat \varphi(z)$  be defined as in
(\ref{non:reg:est:classical}), where
$K$ is a multiplicative kernel of order $r$. In addition assume that 
$\varphi\cdot f,\,f\in \mf{G}^{d}_{s,\triangle}$ 
for some $\triangle,s>0$, $f$ is bounded away from zero and let $\rho=r\wedge s$.
Then we have for all $z\in\R^d$
\begin{equation*}
|\hat \varphi(z)- \varphi(z)|^{2}=
O_{P}\Bigl(b^{2\rho} + b^{-d}\cdot T^{-(\gu\wedge1)} + 
b^{-d(1+q+q_{FG}(1-[(q\gv)\vee 1])}\cdot T^{-[(q\gv)\wedge 1]}  \Bigr), 
\quad T\to\infty.
\end{equation*}
\end{theorem}

\begin{remark}
We observe that the bound obtained in Theorem \ref{non:reg:theo-model1} are 
similar to the bound derived for the density estimator in Theorem
\ref{non:theo}, where
\begin{eqnarray}\label{eq:remark-nonparametrics}
\Ex|\hat f(z)- f(z)|^{2}=O\Bigl(b^{2\rho} + b^{-d}\cdot T^{-1} + 
b^{-d(2-q_{F}(1-[\gv\vee 1]))}\cdot T^{-[\gv\wedge 1]}  \Bigr),
\end{eqnarray}
noting that the result above is for arbitrary dimension $d$. 
The difference is the inclusion of the covariance size $\gu$
of the errors and the $q$ which `balances' the tails of $1/\varphi$ and $f$ 
(see Assumption \ref{non:reg:ass:tec-part1}(i)). 
However, 
we observe that we can partition the bound in Theorem \ref{non:reg:theo-model1} into three cases, which are similar
to the three cases considered in Remark \ref{remark:density}. 
Most notably, we observe
if $\gu>1$ and $\mathfrak{v} >1/q_{FG}+1/q$ then we obtain the usual 
bound $O(b^{2\rho} + b^{-d}\cdot T^{-1})$ for the MSE.

It is interesting to note that in the case
$h^{p}\cdot f$ and $|\varphi|^{p}\cdot f$ are uniformly bounded for all $p$, 
then $q = 1$ (eg. $h$ and $\varphi$ are bounded functions and $f$ is exponential 
density). 
In this case the bounds given in (\ref{eq:remark-nonparametrics}) and 
Theorem \ref{non:reg:theo-model1} are quite similar. The main difference
is the appearance of $q_{FG}$ rather than $q_{F}$ and, the term 
$b^{-d}T^{-(\gu\wedge1)}$ which replaces $b^{-d}T^{-1}$.\hfill $\Box$
\end{remark}

\begin{corollary}\label{non:reg:cor-model1}
Suppose the assumptions of Theorem~\ref{non:reg:theo-model1} 
are satisfied. Let $b^*\approx T^{-\gamma / (2\rho+d)}$ with
\begin{equation}\label{non:reg:cor-model1:gamma}
\gamma:=\left\{\begin{array}{ll}
\min(\gu, 1),&q\gv>1+1/q_{FG};\\
\min\Bigl(\gu, [(q\gv)\wedge 1]
\cdot\frac{2\rho +d}{2\rho+d(1+q+q_{FG}(1-[(q\gv)\vee 1]))}\Bigr),
&1+1/q_{FG}\geq q\gv.
\end{array}\right.,
 \end{equation}
then we have $|\hat \varphi(z)- \varphi(z)|^{2}=O_{P}
\Bigl(T^{-\frac{2\rho}{2\rho+d}\cdot\gamma}\Bigr)$ for all $z\in\R$.
\end{corollary}

Let us now compare Theorem
\ref{non:reg:theo-model1} with the bound 
obtained for the deterministic design $X_{t}= \varphi(\frac{t}{T}) + 
\varepsilon_{t}$, where 
$\gu$ is the covariance size of the errors. In the case of the fixed design, 
the bound for the deviation of the kernel estimator is 
$O(b^{2\rho} + T^{-(\gu\wedge1)}b^{-d})$ (c.f. \cite{p:hal-har-90b}). 
We see that the bound in Theorem \ref{non:reg:theo-model1} include this term,
but also the additional term 
$O(b^{-d(1+q+q_{FG}(1-[(q\gv)\vee 1])}\cdot T^{-[(q\gv)\wedge 1]})$, which
is the influence of the design, in particular, $\gv$. 
If the mixing size of the design were
sufficiently large, then the fixed design and random design estimators have the 
same rate of convergence, $O(T^{-\frac{2\rho}{2\rho+d}})$. 

\begin{remark}[Example]\label{remark:example}
Examples of processes which satisfy (\ref{non:reg:model1})
are stochastic volatility models (c.f  \cite{p:mam-04}), where one observes 
$\{Y_{t}\}$, which satisfies the representation
\begin{eqnarray*}
Y_{t} = \sigma(Z_{t})\eta_{t}.
\end{eqnarray*}
Here $\{\eta_{t}\}$ are iid random variables, 
$\Ex(\eta_{t}^{2})=1$ and $\{Z_{t}\}$ are explanatory variables which can include past values
of $Y_{t}$. Usually in finance the object is to estimate the conditional volatility 
$\sigma^{2}$. By noting that $Y_{t}^{2}$ can be written as
\begin{eqnarray*}
Y_{t}^{2} = \sigma(Z_{t})^{2} + 
(\eta_{t}^{2}-1)\sigma(Z_{t})^{2},  
\end{eqnarray*}
we see that $Y_{t}^{2}$ satisfies  (\ref{non:reg:model1}) with
$X_{t}=Y_{t}^{2}$, 
$\varepsilon_{t} =  (\eta_{t}^{2}-1)$ and $h(\cdot) =\sigma(\cdot)^{2}$. 
Therefore we can estimate the volatility $\sigma(\cdot)^{2}$ using 
(\ref{non:reg:est:classical}), where 
$\hat{\sigma}(\cdot)^{2}$, is the kernel estimator of $\sigma(\cdot)^{2}$.
Furthermore, Theorem \ref{non:reg:theo-model1} can be applied to obtain the rate of
convergence. More precisely, let $\gv$ be the mixing size of $\{Z_{t}\}$, and
noting that 
$\cov\{(\eta_{t}^{2}-1),(\eta_{s}^{2}-1)\} = 0$, when $t\neq s$,
which implies $\gu = \infty$, we obtain 
\begin{equation*}
|\hat \sigma(z)^{2}- \sigma(z)^{2}|^{2}=
O_{P}\Bigl(b^{2\rho} + 
b^{-d(1+q+q_{FG}(1-[(q\gv)\vee 1])}\cdot T^{-[(q\gv)\wedge 1]}  \Bigr).
\qquad \qquad \qquad
\square
\end{equation*}
\end{remark}
From Corollary \ref{non:reg:cor-model1} we see that there are two factors 
which affect the rate of convergence:
the mixing size $\gv$ of the random design $\{Z_{t}\}$ and the size $\gu$
of the covariance
function of $\{\eta_{t}\}$. There are however several models of interest, 
which do not satisfy condition (\ref{non:reg:model1}).
In this case Theorem \ref{non:reg:theo-model1} cannot be applied and it 
is of interest to investigate what happens in the general case. 

Examples of models which do not necessarily
satisfy (\ref{non:reg:model1}) include the 
Cheng-Robinson model, where $\{X_{t}\}$ satisfies the representation
$X_{t} = F(U_{t}) + G(U_{t},Y_{t})$ with $\Ex(G(U_{t},Y_{t})|U_{t})=0$ and   
$\{Y_{t}\}$ is a long memory process, which is 
independent of the weakly dependent design random variables
$\{U_{t}\}$ (c.f. \cite{p:che-rob-94},  
Cs\"org\"o and  Mielniczuk (1999, 2001)\nocite{p:cso-mie-99}\nocite{p:cso-mie-01}).
However, the results are derived under the assumption that $\{Y_{t}\}$ comes from a 
linear process  and $G(\cdot)$ has a particular form.

An alternative approach is developed in \cite{b:bos-98}, who considers nonparametric 
prediction for time series, where one observes the stationary time series
$\{(X_{t},Z_{t})\}$ and the parameter of interest is $\varphi(z) = 
\Ex(X_{t}|Z_{t}=z)$. The sampling results in \cite{b:bos-98}
are based on the assumption that the mixing size of $\{(X_{t},Z_{t})\}$
is sufficiently large, (thus excluding Cheng-Robinson type models) yielding 
an estimate which has the same rate as the kernel estimator for 
iid random variables.

We now consider the sampling properties of $\hat{\varphi}$, when the observations 
$\{(X_{t},Z_{t})\}$ satisfy the general model defined in (\ref{non:reg:model}), 
and dependence is quantified through its 2-mixing size, which can be arbitrary. 

We will use the following assumptions.
\begin{assumption}[Densities, moments and kernels]
\label{non:reg:ass:tec-part2}\begin{itemize}
\item[(i)] Let $\ex|X_{t}|^p<\infty$ for some $p>2$ and define  
$g^{(p)}(z):=\ex[|X_{t}|^p\vert Z_{t}=z]\cdot f(z)$. Then the functions $g^{(p)}$ and $f$ are 
uniformly  bounded and we define  $q:=1-2/p$. 
\item[(ii)]Let 
$f^{(t,\tau)}$ and $F^{(t,\tau)}$ be defined as in 
Assumption \ref{dens:est:ass:tec-part} (ii) and let 
$g^{(t,\tau)}$ and $G^{(t,\tau)}$ be defined as in 
Assumption \ref{dens:est:ass:tec-part} (ii).
Then $\|F^{(t,\tau)}\|_{p_{F}}$ and $\|G^{(t,\tau)}\|_{p_{G}}$ are uniformly bounded in $t$ and $\tau$ for some $p_{F},p_{G}>2$, where we define $q_{F}:=1-2/p_{F}$, $q_{G}:=1-2/p_{G}$ and $q_{FG}:=q_{F}\wedge q_{G}$.
\item[(iii)] The multiplicative kernel $K$ has finite first and $p$-th moment.
\end{itemize}
\end{assumption}
We note that assumptions above are similar to
Assumption \ref{non:reg:ass:tec-part1}. The difference lies in 
Assumption \ref{non:reg:ass:tec-part1}(i) and 
Assumption \ref{non:reg:ass:tec-part2}(i).
Assumption \ref{non:reg:ass:tec-part2}(i) is in terms of moments
whereas Assumption \ref{non:reg:ass:tec-part1}(i) is in terms of functions.

In the following theorem we derive an error bound  for the estimator 
$\hat\varphi$.
\begin{theorem}\label{non:reg:theo-model2}
Suppose the stationary time series $\{(X_{t},Z_{t})\}$ satisfies 
(\ref{non:reg:model}), and is  
$2$-mixing  of size $\gv$.
Furthermore, Assumption  \ref{non:reg:ass:tec-part2} is fulfilled for some 
$ q_{FG},q\in(0,1)$. 

Let the estimator $\hat \varphi(z)$ be defined as in 
(\ref{non:reg:est:classical}), where
$K$ is a multiplicative kernel of order $r$. In addition assume that 
$\varphi\cdot f,\,f\in \mf{G}^{d}_{s,\triangle_{}}$ 
for some $\triangle,s>0$,  
$f$ is bounded away from zero and let $\rho= r\wedge s$. Then we
have for all $z\in\R^d$
\begin{equation*}
|\hat \varphi(z)- \varphi(z)|^{2}=O_{P}\Bigl(b^{2\rho} + 
b^{-d}\cdot T^{-1} + b^{-d(1+q+q_{FG}(1-[(q\gv)\vee 1])}\cdot
T^{-[(q\gv)\wedge 1]}  
\Bigr), \quad T\to\infty.
\end{equation*}
\end{theorem}
We now obtain the rates of convergence using the optimal bandwidth. 
\begin{corollary}\label{non:reg:cor-model2}
Suppose the assumptions in Theorem~\ref{non:reg:theo-model2} are satisfied.  
Let  $b^*\approx T^{-\gamma / (2\rho+d)}$ with
\begin{equation}\label{non:reg:cor-model2:gamma}
\gamma:=\left\{\begin{array}{ll}
1,&q\gv>1+q/q_{FG};\\
\,[(q\gv)\wedge 1]\cdot\frac{2\rho +d}{2\rho+d(1+q+q_{FG}(1-[(q\gv)\vee 1]))}
,&1+q/q_{FG}\geq q\gv., 
\end{array}\right.
 \end{equation}
 then we have $|\hat \varphi(z)- \varphi(z)|^{2}
=O_{P}\Bigl(T^{-\frac{2\rho}{2\rho+d}\cdot\gamma}\Bigr)$ for all $z\in\R$.
\end{corollary}

%% file: panel.tex
\section{Nonparametric panel time series}\label{sec:non:pan}
In recent years, panel time series have often been used to model the relationship
and dynamics between several observed time series. 
Typically we let $X_{t,i}$ denote the observation of the $i$th  
individual at time $t$, where $i=1,\ldots,N$ and $t=1,\ldots,T$. 
 We also assume that we observe some explanatory variables $Z_{t,i}$
which are known to influence $X_{t,i}$.
Several models have been proposed 
to model the complex relationship between individuals,
ranging from parametric models (c.f  \cite{b:bal-01}, \cite{p:tjo-99}, 
\cite{p:dah-fei-05}, 
and the references therein) to nonparametric additive models
(c.f. \cite{p:mam-05}). 
In this section, we take the nonparametric route, and use the methods developed in 
the sections above to obtain an estimator of the mean function and study its
sampling properties.  The results in this section can be used in various
applications, an interesting example is the estimation of the covariance 
function of spatial-temporal models considered in \cite{p:joh-jun-sub-07}.  

Let us suppose the affect of the explanatory 
variables is common over all individuals. To be precise, the response and 
explanatory variables $\{(X_{t,i},Z_{t,i})\}_{t}$ form a ($1+d$)-dimensional 
stationary vector time series which satisfies the relation
\begin{equation}\label{non:pan:model}
\ex[X_{t,i}|Z_{t,i}=z]=\varphi(z)\quad\forall z\in\R^d, i\in\N, t\in \Z.
\end{equation}
We describe the dependence of $\{(X_{t,i},Z_{t,i})\}$, by assuming it is 
$2$-mixing over time, see Definition \ref{non:pan:def:t-dep}, below. 
We note that the model considered in \cite{p:tjo-04} and \cite{p:mam-05}, 
can be used as a particular example of (\ref{non:pan:model}).

We now define an estimator for $\varphi$. Note that we do not suppose that  
different individuals, say $(X_{t,i},Z_{t,i})$ and $(X_{t,j},Z_{t,j})$
are identically distributed  (have common densities). 
Let $f_{i}(x,z)$ denote the joint density of the random vector 
$(X_{t,i},Z_{t,i})$ for $i\in \N$.  
Moreover, let  $f_{i}(z)$  denote the marginal density of  $Z_{t,i}$. 
Using these densities we can rewrite (\ref{non:pan:model}) as
\begin{displaymath} \varphi(z)=\ex[X_{t,i}|Z_{t,i}=z]=
\int x \frac{f_{i}(x,z)}{f_{i}(z)}dx=
:\frac{g_{i}(z)}{f_{i}(z)},\quad \forall z\in\R^d,t\in\Z,i\in\N.
\end{displaymath}
Furthermore, using the above,  it is easily verified that
\begin{equation}\label{non:pan:model:1}
\varphi(z)=\frac{\frac{1}{N}
\sum_{i=1}^Ng_{i}(z)}{\frac{1}{N}\sum_{i=1}^Nf_{i}(z)}, 
\quad \forall z\in\R^d,N\in\N
\end{equation}
which motivates the following estimator of $\varphi$. 

Given the observations 
$\{(X_{t,i},Z_{t,i});t=1,\dotsc,T;i=1,\dotsc,N\}$ our object is to estimate 
$\varphi$ and consider its sampling properties. 
The identity (\ref{non:pan:model:1}) suggests as an estimator of 
$\varphi(z)$ 
\begin{equation}\label{non:pan:est}
\hat\varphi(z)=\frac{\frac{1}{N}\sum_{i=1}^N
\hat g_{i}(z)}{\frac{1}{N}\sum_{i=1}^N\hat f_{i}(z)},
\end{equation}
using  for each $i=1,\dotsc,N$
\begin{gather}\label{non:pan:est:gf}
\hat g_{i}(z):=\frac{1}{T}\sum_{t=1}^TX_{t,i}K_{b_{i}}(Z_{t,i}-z)
\quad\mbox{and}\quad
\hat f_{i}(z):=\frac{1}{T}\sum_{t=1}^TK_{b_{i}}(Z_{t,i}-z)
\end{gather}
as estimators of $g_{i}(z)$ and $f_{i}(z)$, respectively. 

We quantify the dependence both over time and between individuals through
their 2-mixing rates.
\begin{defin}\label{non:pan:def:t-dep} The panel time series 
$\{(X_{t,i},Z_{t,i})\}_{t}$, 
$i\in\N$, is said to be $2$-mixing with size $\gv $ and $\gu$, if for all $i,j\in \N$
\begin{displaymath}
\sup_{A\in \sigma(X_{t,i},Z_{t,i}),B\in
\sigma(X_{\tau,j},Z_{\tau,j})} |P(A\cap B) - P(A)P(B)|\leq
C \left\{\begin{array}{ll}|t-\tau|^{-\mathfrak{v}};&
\mbox{if }i=j,\\|t-\tau|^{-\mathfrak{u}};&\mbox{otherwise},\end{array}\right.
\end{displaymath}
for some $C<\infty$ independent of $i,$ $j$, $t$ and $\tau$.
\end{defin}
In the results below we will show that the rate of convergence of $\varphi(\cdot)$
is determined by the smallest mixing size $(\gv\wedge \gu)$. However in the 
case that $\gv < \gu$ and the number of individuals $N$ grow with 
$T$, the rate is determined, solely, by $\gu$. 

\begin{remark}\label{remark:pan:individuals}
We note that in Definition \ref{non:pan:def:t-dep} we have two 
2-mixing sizes, the size $\mathfrak{v}$ describes the dependence of the 
time series $\{(X_{t,i},Z_{t,i})\}_{t}$, whereas the size $\mathfrak{u}$ describes the
dependence between individuals over time. By separating these two sizes we
can model different behaviours. A simple example is  $X_{t,i} = \varphi(Z_{i}) + \varepsilon_{t,i}$, where for a given 
individual $i$, the 
explanatory variable $Z_{i}$ is fixed over time, 
 $\{\varepsilon_{t,i}\}$ and $\{Z_{i}\}$ are iid random variables. In this example,  
$\mathfrak{v} = 0$ and $\mathfrak{u}=\infty$.  \hfill $\Box$
\end{remark}

To obtain the sampling properties of $\hat \varphi$ we require the following assumptions,
which are an extension of Assumption \ref{non:reg:ass:tec-part2} to panel data. 
\begin{assumption}[Densities, moments and kernels]\label{non:pan:ass:tec}
\begin{itemize}
\item[(i)]For all $i\in\N$ let $\ex[|X_{t,i}|^{p}]<\infty$ for some $p>2$ and define 
$g_{i}^{(p)}(z) := \ex[X_{t,i}^p| 
Z_{t,i}=z ]\cdot f_{i}(z)$. Then the functions $g_{i}^{(p)}$ and $f_{i}$, for
all $i\in \N$, are uniformly bounded  and we define  $q:=1-2/p$. 
\item[(ii)]For each $t,\tau\in\Z$ and $i,j\in\N$ let 
$f_{i,j}^{(t,\tau)}$ denote
the joint density of $(Z_{t,i},Z_{\tau,j})$  and let
$g^{(t,\tau)}_{i,j}(z_{1},z_{2}):= \ex[X_{t,i}X_{\tau,j}| Z_{t,i}=z_{1},
Z_{\tau,j}=z_{2} ]\cdot f_{i,j}^{(t,\tau)}(z_{1},z_{2})$. Define
$F_{i,j}^{(t,\tau)} :=  f_{i,j}^{(t,\tau)} - f_{i} \otimes f_{j}$ and   
$G_{i,j}^{(t,\tau)} := g_{i,j}^{(t,\tau)} - g_{i} \otimes g_{j}$. Then 
$\|F_{i,j}^{(t,\tau)}\|_{p_{F}}$ and $\|G_{i,j}^{(t,\tau)}\|_{p_{G}}$ are uniformly 
bounded in $i,j,t$ and $\tau$ for some $p_{F},p_{G}>2$. We define $q_{F}=1-2/p_{F}$, 
$q_{G}=1-2/p_{G}$ and $q_{FG}:=q_{F}\wedge q_{G}$.
\item[(iii)] The multiplicative kernel $K$ has a finite first and $p$-th moment.
\end{itemize}
\end{assumption}
We now obtain a bound for the deviation of $\hat \varphi(z)$, as
$N$ is kept fixed and $T\rightarrow \infty$. 

\begin{theorem}\label{non:pan:theo:MSE} Let us suppose that 
the stationary panel time series
$\{(X_{t,i},Z_{t,i})\}$ satisfies (\ref{non:pan:model}),  
for all $i,j\in \N$, and is 
$2$-mixing with size $\mathfrak{v}$ and $\mathfrak{u}$ 
(as defined in Definition \ref{non:pan:def:t-dep}). 
Suppose Assumption \ref{non:pan:ass:tec} is fulfilled for some $ q,q_{FG}\in(0,1)$. 

Let the estimator $\hat \varphi$ be defined in (\ref{non:pan:est}), where  
for each $i=1,\dotsc,N$ 
the nonparametric estimators $\hat g_{i}$ and $\hat f_{i}$ given in 
(\ref{non:pan:est:gf}) are  
constructed using a  multiplicative kernel $K$ of order $r>0$. 
In addition assume for each $i=1,\dotsc,N$, that the functions 
$\varphi\cdot f_{i}$ and $f_{i}$ 
belong to $ \mf{G}^{d}_{s_{i},\triangle}$ for $s_{i},\triangle>0$,  
$f_{i}$ is bounded away from zero and let $\rho_{i}=r\wedge s_{i}$. 
Then we have for all $z\in\R^{d}$
\begin{multline}
\label{non:pan:theo:MSE:e}
 |\hat \varphi(z)- \varphi(z)|^{2} = 
O_{p} \Bigl(\frac{1}{N}
\sum_{i=1}^N\Bigl\{ b_{i}^{2\rho_{i}} +T^{-1}\cdot b_{i}^{-d}+ T^{-[(q\gu)\wedge 1]} \cdot 
 b_{i}^{-d(1+q+q_{FG}-q_{FG}[(q\mathfrak{u})\vee 1])}\\
\hfill+
N^{-1}\cdot T^{-[(q\gv)\wedge 1]} \cdot 
 b_{i}^{-d(1+q+q_{FG}-q_{FG}[(q\mathfrak{v})\vee 1])}\Bigr\}\Bigr),\quad T\to\infty.
\end{multline}
\end{theorem}
Comparing Theorem \ref{non:reg:theo-model2} with the theorem above, we see, 
besides the  summation $\frac{1}{N}\sum_{i=1}^N$,  the addition of an extra term 
$T^{-[(q\gu)\wedge 1]} \cdot b_{i}^{-d(1+q+q_{FG}-q_{FG}[(q\mathfrak{u})\vee 1])}$ 
due to the dependence between individuals over time. Altogether this implies 
that the bound 
can be partitioned into nine different cases, depending on the values of 
$\gu$ and $\gv$
(compare this with the bound in Theorem
\ref{non:reg:theo-model2}, which can be partitioned into three cases). 
However,  the nine different bounds can be grouped into
two main cases; when $\gu \leq \gv$ and $\gu > \gv$, we consider these two 
cases in the corollaries below. 

If $\gu \leq \gv$, we notice that the third term dominates the fourth term, 
in other words there is 
a larger dependence between individuals over time than for each individual
over time. We consider this case below. 

\begin{corollary}\label{non:pan:cor:MSE:gu} 
Suppose the assumptions in Theorem~\ref{non:pan:theo:MSE}
 are satisfied and $\gu\leq \gv$. For each $i=1,\dotsc,N$, let  
$b_{i}^*\approx T^{-\gamma_{i} / (2\rho_{i}+d)}$ with
\begin{equation}\label{non:pan:cor:MSE:gu:gamma}
\gamma_{i}:=\left\{\begin{array}{ll}
1,&q\gu>1+q/q_{FG};\\
\,[q\gu\wedge1]\cdot\frac{2\rho_{i} +d}{2\rho_{i}+d(1+q+q_{FG}(1-[(q\gu)\vee1]))},&1+q/q_{FG}\geq q\gu,
\end{array}\right. 
 \end{equation}
then for all $z\in\R^d$ we have   
$|\hat \varphi(z)- \varphi(z)|^{2}=O_{P}\Bigl(\frac{1}{N}
\sum_{i=1}^N T^{-\frac{2\rho_{i}}{2\rho_{i}+d}\cdot\gamma_{i}}\Bigr),$ $T\to\infty$.
\end{corollary}
Studying the corollary above we see that the rate of convergence is 
determined by $\gu$. 
In other words, there is no benefit in the estimation by including several 
individuals $N$. Furthermore we see that the `usual' nonparametric rate 
of convergence is only achieved if $\gu>1/q+1_{FG}$. 

Let us now consider the situation where $\gv< \gu$, that is there is less 
dependence between two different individuals 
over time than one individual observed over time (see Remark
\ref{remark:pan:individuals} for an example). This scenario 
is more likely to arise in real applications. In this case, the fourth term 
dominates the third term in (\ref{non:pan:theo:MSE:e}), which suggests that
increasing the number of individuals
does yield a faster rate of convergence.   We notice that the
usual nonparametric rate of convergence can only be obtained if 
$\gu > 1/q_{FG} +1/q$. 

\begin{corollary}\label{non:pan:cor:MSE:gr} Suppose the assumptions of 
Theorem~\ref{non:pan:theo:MSE} 
are satisfied and $\gv\leq \gu$. For each $i=1,\dotsc,N$, let  
$b_{i}^*\approx N^{-\zeta_{i} / (2\rho_{i}+d)}\cdot T^{-\delta _{i} / (2\rho_{i}+d)}$ 
where $\zeta_{i}\geq 0$ and 
\begin{equation}\label{non:pan:cor:MSE:gr:delta}
\delta_{i}:=\left\{\begin{array}{ll}
1,&q\gv>1+q/q_{FG};\\
\,[q\gv\wedge1]\frac{2\rho_{i} +d}{2\rho_{i}+d(1+q+q_{FG}(1-[q\gv\vee1]))},&
1+q/q_{FG}\geq q\gv.
\end{array}\right.
 \end{equation}
 Then for all $z\in\R^d$ we have  
 \begin{multline}\label{non:pan:cor:MSE:gr:e}
  |\hat \varphi(z)- \varphi(z)|^{2}=O_{P}\Bigl(\frac{1}{N}
\sum_{i=1}^N\Bigl\{ N^{-\frac{2\rho_{i}}{2\rho_{i} + d}\cdot \zeta_{i}}\cdot
T^{-\frac{2\rho_{i}}{2\rho_{i}+d}\cdot\delta_{i}}\cdot 
[N^{\frac{\zeta_{i}}{\gamma_{i}}\cdot[(q\gu)\wedge1]}\cdot T^{-\frac{\gamma_{i}-\delta_{i}}{\delta_{i}}\cdot[(q\gu)\wedge1]}+ \\
\hfill N^{-1+\frac{\zeta_{i}}{\delta_{i}}\cdot[(q\gv)\wedge1]}]\Bigr\}
\Bigr),
\quad T\to\infty, \end{multline}
where $\gamma_{i}$ is defined in (\ref{non:pan:cor:MSE:gu:gamma}).
\end{corollary}

Studying the corollary above, we see if $N$ is kept fixed, then the terms
inside the inner bracket of (\ref{non:pan:cor:MSE:gr:e})
are of order $O_{P}(1)$, therefore the rate of convergence is $O_{P}\Bigl(\frac{1}{N}
\sum_{i=1}^N T^{-\frac{2\rho_{i}}{2\rho_{i}+d}\cdot\delta_{i}}\Bigr)$. Thus, 
combining Corollaries  \ref{non:pan:cor:MSE:gu} and \ref{non:pan:cor:MSE:gr}
we have for arbitrary $\gv $ and 
$\gu$ the rate\linebreak $O_{P}\Bigl(\frac{1}{N}
\sum_{i=1}^N
T^{-\frac{2\rho_{i}}{2\rho_{i}+d}\cdot[\gamma_{i}\wedge\delta_{i}]}
\Bigr)$. Therefore, we obtain the usual nonparametric rate if 
$(\gv\wedge \gu)> 1/q_{FG} +1/q$.

Altogether the corollaries above imply that the rate of convergence depends on 
the slowest mixing rate, within or between the individuals. Let us suppose
that $\gv < \gu$, 
if we closely examine (\ref{non:pan:cor:MSE:gr:e}) we see if $\zeta_{i}$ is
chosen such that
$\zeta_{i}<\gamma_{i}[(q\gu)\wedge 1)]$, $\gamma_{i}<\delta_{i}[(q\gu)\wedge 1]$ 
(noting that the former inequality implies the later, since $\gv < \gu$)
and $(\delta_{i} - \gamma_{i})>0$ (which is the case when 
$\gv < \gu$), then the terms inside the inner bracket of
(\ref{non:pan:cor:MSE:gr:e}) become small for large $N$. 
This means that increasing the number of individuals  leads to a faster rate of 
convergence. We show in the corollary below if we allow the number of individuals 
$N$ to grow as $T$ grows,   
then the rate of convergence will depend only 
on $\gu$ and no longer on the smaller $\gv$ (unlike the case that  $N$ is  fixed). 

\begin{corollary}\label{non:pan:cor:MSE:gr:N} Suppose the assumptions in 
Theorem~\ref{non:pan:theo:MSE} 
are satisfied and let $\gv\leq \gu$. Furthermore assume there exists a 
$\zeta_{i} > 0$  such that for each $i\in\N$, 
$N^{\zeta_{i}}\approx T^{(\gamma_{i}-\delta_{i})}$, where $\gamma_{i}$ and 
$\delta_{i}$ are defined in 
(\ref{non:pan:cor:MSE:gu:gamma}) and (\ref{non:pan:cor:MSE:gr:delta}),
respectively and  $\delta_{i}/[(q\gv)\wedge 1]\geq \zeta_{i}$. 
Then given 
$b_{i}^*\approx N^{-\zeta_{i}/(2\rho_{i}+d)}\cdot T^{-\delta_{i}/(2\rho_{i}+d)} 
= T^{-\gamma_{i} /(2\rho_{i}+d)}$  we have  for all $z\in\R^d$
 \begin{equation}\label{non:pan:cor:MSE:gr:N:e}
  |\hat \varphi(z)- \varphi(z)|^{2}=O_{P}\Bigl(\frac{1}{N}
\sum_{i=1}^N  T^{-\frac{2\rho_{i}}{2\rho_{i}+d}\cdot\gamma_{i}}\Bigr),\quad T\to\infty. \end{equation}
\end{corollary}
We see from the corollary above if the number of individuals, $N$, grows at the rate $N 
\approx T^{\frac{\gamma_{i}-\delta_{i}}{\zeta_{i}}}$, where $\zeta_{i}$ 
cannot be too large, in particular
$\delta_{i}[(q\gv)\wedge 1]\geq \zeta_{i}$, then the rate of convergence 
depends only on the mixing size
$\gu$ (compare this with Corollaries \ref{non:pan:cor:MSE:gu} 
and \ref{non:pan:cor:MSE:gr}, 
where the mixing size depends on $(\gv \wedge \gu)$). 
Furthermore if $q\gu\geq q/q_{FG}+1$ then 
we have the usual nonparametric rate 
$|\hat \varphi(z)- \varphi(z)|^{2}=O_{P}\Bigl(\frac{1}{N}
\sum_{i=1}^N  T^{-\frac{2\rho_{i}}{2\rho_{i}+d}}\Bigr)$. 
In the special case that  $\varphi$ and all the marginal densities $f_{i}$ 
belong to the same smoothness class,
that is for all $i$, $\rho =\rho_{i}$,  then (\ref{non:pan:cor:MSE:gr:N:e}) 
simplifies to 
$|\hat \varphi(z)- \varphi(z)|^{2}=
O_{P}\Bigl( T^{-\frac{2\rho}{2\rho+d}}\Bigr)$. 

\begin{remark}
It is worth mentioning that similar results to those in Theorem  
\ref{non:pan:theo:MSE} and Corollaries \ref{non:pan:cor:MSE:gu}, 
\ref{non:pan:cor:MSE:gr} and \ref{non:pan:cor:MSE:gr:N} can be obtained for
the model (\ref{non:reg:model1}) considered in Section \ref{sec:non:reg}. 
\hfill $\Box$
\end{remark}

%% file: discussion.tex
\section{Discussion}
In this paper we have considered nonparametric estimation for 
dependent data. Focusing on the case that the observations are 
nonlinear and highly dependent. We have obtained bounds for the 
kernel density estimator and also rates of convergence
of two types of nonparametric regression models, both using the 
2-mixing dependence measure.  We show that when the assumption 
of linearity is relaxed, the rate of convergence does not necessarily
depend on the autocovariance function of the observations. 
We demonstrate that 2-mixing 
is a natural measure of dependence for panel data and obtained
rates of convergence for the common mean function in panel time 
series. 

As we are working under relatively weak conditions, 
we do not claim that the bounds obtained are minimax. However, the bounds can be 
considered as the worst case scenario for the nonparametric estimator. In 
future work, it would be of interest to investigate if the bounds in the paper are 
indeed close to minimax for certain nonlinear time series.  
In this paper we have derived
bounds for the estimator using the optimal bandwidth. However the optimal
bandwidth is constructed under the assumption that the 2-mixing size is known.  
It would also be on interest to develop bandwidth selection methods 
when the 2-mixing size of the observations is unknown. 

%% file: appendix-proofs-nonpar-density.tex
\subsection{Proofs: Nonparametric density estimation}\label{a:p:non:dens}
We now prove the results in Section \ref{sec:density}. We mention that 
Theorem \ref{non:theo}
is stated for an univariate time series $\{Z_{t}\}$, however the 
proofs of the results in Sections \ref{sec:non:reg} and \ref{sec:non:pan} 
require results in the multivariate case. 
Therefore to save space, we give the proof of
Theorem \ref{non:theo} for a $d$-dimensional vector time series $\{Z_{t}\}$. 

\begin{lemma}\label{a:p:non:dens:l:1} 
Suppose the time series $\{Z_{t}\}$ is $2$-mixing with size  $\mathfrak{v}$ and 
 Assumption \ref{dens:est:ass:tec-part} is fulfilled for some $q_{F}\in(0,1)$. 
If $1\leq t,\tau\leq T$, then
\footnote{We write $A \lesssim B$ is there exists a positive constant $c$ such that $A \leqslant c B$.}
\begin{multline}\label{a:p:non:dens:l:1:res1}
|\cov\left\{K_{b}(Z_{t}-z),K_{b}(Z_{\tau}-z)\right\}|  
\lesssim\min\Bigl(b^{-d(1-q_{F})}; b^{-2d}|t-\tau|^{-\mathfrak{v}}
\Bigr).
\end{multline}
\end{lemma}

\bigskip

\noindent\textcolor{darkred}{\sc Proof.} Writing the covariance as an integral, and using the notation in 
Assumption \ref{dens:est:ass:tec-part} (ii) we have
\begin{displaymath}
 \cov\left\{K_{b}(Z_{t}-z),
K_{b}(Z_{\tau}-z)\right\} 
= \int\; K_{b}(u-z)
K_{b}(v-z)F^{(t,\tau)}(u,v)dudv.
\end{displaymath}
Now by using H\"older's inequality with $p_{F}^{-1}+\bar p_{F}^{-1}=1$, and recalling that
$q_{F} = 1 - 2/p_{F}$,  it is clear that 
\begin{displaymath}
|\cov\left\{ K_{b}(Z_{t}-z),
K_{b}(Z_{\tau}-z)\right\} | 
\leq 
\frac{1}{b^{2d}}\cdot
b^{2d/\bar{p}_{F}} \|K\|_{\bar p_{F}}^2\cdot \|F^{t,\tau}\|_{p_{F}} \lesssim b^{-d(1-q_{F})}.
\end{displaymath}
Using Assumption \ref{dens:est:ass:tec-part} we have that 
 $\|F^{t,\tau}\|_{p_{F}}$ is uniformly bounded and by using Lyaponov's inequality
$\|K\|_{\bar{p}_{F}}<\infty$  for all $1<{\bar p_{F}}<2$. This gives us the  first bound 
in  (\ref{a:p:non:dens:l:1:res1}). On the other hand, under 
Assumption \ref{dens:est:ass:tec-part} (i) the kernel $K$ is uniformly bounded and therefore, 
using the $2$-mixing property of 
$\{Z_{t}\}$ together with  \cite{b:hal-hey-80}, Theorem A.5, we obtain
\begin{displaymath}
|\cov\left\{ K_{b}(Z_{t}-z),K_{b}(Z_{\tau}-z)\right\} | \lesssim b^{-2d} \cdot |t-\tau|^{-\mathfrak{v}},
\end{displaymath}
which gives the second bound in (\ref{a:p:non:dens:l:1:res1}).\hfill $\square$

\bigskip

\noindent\textcolor{darkred}{\sc Proof of Theorem~\ref{non:theo}.}
We mention that parts of the following
proof are motivated by techniques
used in \cite{b:bos-98}, where nonparametric smoothing was considered for
 univariate time series. Consider the standard variance bias decomposition\begin{equation}
 \ex|\hat f(z)- f(z)|^2 = \var(\hat f(z)) + 
|\ex\hat f(z)-f(z)|^2.
\end{equation}
Under the stated assumptions we will derive the following two bounds, which give together the result of the theorem.  The bias is bounded by
\begin{gather}\label{a:p:non:dens:th:MSE:bias}
|\ex\hat f(z)-f(z)|^2\lesssim  b^{2\rho},
\end{gather}
while for the variance we have
\begin{equation}\label{a:p:non:dens:th:MSE:var}
\var(\hat f(z))\lesssim T^{-1} \cdot   b^{-d}+  T^{-[\gv \wedge 1]} \cdot b^{-d(2+q_{F}-q_{F} [\mathfrak{v}\vee 1])}.
\end{equation}

Proof of (\ref{a:p:non:dens:th:MSE:bias}). We can write
\begin{displaymath}
\ex\hat f(z)= \frac{1}{T}\sum _{t=1}^{T}\ex\Bigl(K_{b}(Z_{t}-z)\Bigr)=
\int\; du\; f(u) K_{b}(u-z).
\end{displaymath}
Since $f\in \mf{G}^{d}_{s,\triangle}$  and $K$ is a multiplicative kernel  
of order $r$ with $\int du |u|^r K(u)\leq S_{K}$, using a Taylor expansion 
up to the order $\rho=\min(r,s)$ leads to 
$\ex\hat f(z)=f(z)+b^{\rho} R$ with reminder 
$|R|\leq \triangle S_{K}<\infty,$ which proves  
(\ref{a:p:non:dens:th:MSE:bias}).

In order to proof 
(\ref{a:p:non:dens:th:MSE:var}), we consider the expansion
\begin{align}\nonumber
\var(\hat{f}(z)) &= \frac{1}{T^{2}} \sum_{t=1}^{T} 
\var\left\{ K_{b}(Z_{t}-z)\right\}+ \frac{2}{T^{2}} \sum_{t>\tau} 
\cov\left\{K_{b}(Z_{t}-z),K_{b}(Z_{\tau}-z)\right\}\\\label{a:p:non:dens:th:MSE:var:dec}
&=:A_{1}+A_{2}.
\end{align}
We will show that  
$ |A_{1}|\lesssim T^{-1} \cdot  b^{-d}$ and  
\begin{align}\label{a:p:non:dens:th:MSE:var:dec:A2}
|A_{2}|\lesssim\left\{
\begin{array}{ll}
  T^{-\gv} \cdot b^{-2d}, & \mathfrak{v}\leq 1; \\
 T^{-1}\cdot \{   b^{-d} + b^{-d(2+q_{F}-q_{F}\mathfrak{v})}\}, & 
1<\mathfrak{v}.
\end{array}
\right.
\end{align}
Furthermore, if $0 \leq \mathfrak{v} \leq1/q_{F}+1$ then $|A_{1}|$ is dominated by $|A_{2}|$. Whereas for $\mathfrak{v}>1/q_{F}+1$ 
the terms $|A_{1}|$ and $|A_{2}|$ are of the same order $O(T^{-1}b^{-d})$. Therefore, the bounds derived for 
$|A_{2}|$ will lead to  (\ref{a:p:non:dens:th:MSE:var}). 

First let us consider $A_{1}$. Due to  stationarity, 
we have the bound 
\begin{align*}
 T\cdot A_{1} \leq \ex[K_{b}^2(Z_{1}-z)]&= 
\int \;du\; f(u)K_{b}^2(u-z).
\end{align*}  
Since under the stated assumptions $\normV{K}_{2}<\infty$ and the density $f$ is uniformly  bounded 
 this leads to $A_{1}\lesssim T^{-1}\cdot b^{-d}$. 

The term $T\cdot |A_{2}|$ is  bounded by the sum
$4\sum_{t=2}^T |\cov\left\{K_{b}(Z_{t}-z),K_{b}(Z_{1}-z)\right\}|.$
If $\mathfrak{v}\leq 1 $ then we 
estimate the sum  using the second bound in Lemma~\ref{a:p:non:dens:l:1}, i.e., 
$T\cdot |A_{2}|\lesssim T^{-\mathfrak{v} +1} b^{-2d}$, 
which is the first bound in (\ref{a:p:non:dens:th:MSE:var:dec:A2}).  On the other hand if $\mathfrak{v}> 1$ 
we partition the sum into two parts which we estimate separately using  the bounds in  Lemma~\ref{a:p:non:dens:l:1}, thus giving us
\begin{align*}
T\cdot |A_{2}|  \lesssim & \Bigl\{\sum_{t=2}^{h}
 b^{-d(1-q_{F})}   + \sum_{t=h+1}^T b^{-2d} t^{-\mathfrak{v}}\Bigr\}\lesssim \Bigl\{h \cdot  b^{-d(1-q_{F})} +  h^{-\mathfrak{v}+1}\cdot b^{-2d}\Bigr\}.
\end{align*}
Thereby using $h\approx b^{-d q_{F}}$ we obtain  $T \cdot |A_{2}| \lesssim  b^{-d} + b^{-d(2+q_{F}-q_{F}\gv)}$, i.e., the second bound in (\ref{a:p:non:dens:th:MSE:var:dec:A2}). Thus we have proved (\ref{a:p:non:dens:th:MSE:var}).\hfill $\square$

\bigskip

\noindent\textcolor{darkred}{\sc Proof of Corollary~\ref{cor:theo}}
Under the assumption on the bandwidth the result is obtained by balancing the terms in the bound given in Theorem~\ref{non:reg:theo-model1}.\hfill $\square$

%% file: appendix-proofs-nonpar-regression.tex
\subsection{Proofs: Nonparametric regression}\label{a:p:non:reg}
We now prove the results in Section \ref{sec:non:reg}. 
\begin{lemma}\label{a:p:non:reg:l:1} 
Suppose the stationary time series $\{X_{t},Z_{t}\}$ satisfies  (\ref{non:reg:model1}), and $\{Z_{t}\}$ is $2$-mixing with size $\gv$ and the autocovariances of the  
time series $\{\eta_{t}\}$ have size $\gu$ (see Definition \ref{defin:size}).  
Suppose Assumption \ref{non:reg:ass:tec-part1} is fulfilled for some 
$q,q_{G}\in(0,1)$.  If $1\leq t,\tau\leq T$, then
\begin{gather}\label{a:p:non:reg:l:1:res1}
|\cov\left\{X_{t}K_{b}(Z_{t}-z),X_{\tau}K_{b}(Z_{\tau}-z)\right\}|
\lesssim b^{-d(1-q_{G})},\\\label{a:p:non:reg:l:1:res2}
|\cov\left\{\varphi(Z_{t})K_{b}(Z_{t}-z),\varphi(Z_{\tau})K_{b}(Z_{\tau}-z)\right\}|
 \lesssim b^{-d(1+q)} |t-\tau|^{-q\gv},\\\label{a:p:non:reg:l:1:res3}
|\cov\left\{h(Z_{t})K_{b}(Z_{t}-z)\eta_{t},h(Z_{\tau})K_{b}(Z_{\tau}-z)\eta_{\tau}\right\}|
 \lesssim b^{-d} |t-\tau|^{-\gu}.\end{gather}
\end{lemma}

\bigskip

\noindent\textcolor{darkred}{\sc Proof.} Using the notation in 
Assumption \ref{non:reg:ass:tec-part1}  together with H\"older's inequality, and recalling that
$q_{G} = 1 - 2/p_{G}$  with $p_{G}^{-1}+\bar p_{G}^{-1}=1$,  we have 
\begin{displaymath}
|\cov\left\{ X_{t}K_{b}(Z_{t}-z),
X_{\tau}K_{b}(Z_{\tau}-z)\right\} |  \lesssim b^{-d(1-q_{G})},
\end{displaymath}
where we use that under Assumption \ref{non:reg:ass:tec-part1}, 
  $K$ has finite $1< \bar p_{G}< p$ moment (byÊ Lyaponovs inequality) and $\|G_{t,\tau}\|_{p_{G}}$ is uniformly bounded. This gives us  (\ref{a:p:non:reg:l:1:res1}). 

We now prove (\ref{a:p:non:reg:l:1:res2}).  Under Assumption   \ref{non:reg:ass:tec-part1}
the function $|\varphi|^{p}\cdot f$ is uniformly bounded and $\|K\|_{p}$ is finite  
for some $p=2/(1-q)>2$, therefore we have 
$[\ex |\varphi(Z_{1})K_{b}(Z_{1}-z)|^p]^{2/p}\lesssim b^{-d(q+1)}.$
Using the $2$-mixing property of 
$\{Z_{t}\}$ together with  \cite{b:hal-hey-80}, Theorem A.6, we obtain
(\ref{a:p:non:reg:l:1:res2}).

We now
prove (\ref{a:p:non:reg:l:1:res3}). The series $\{Z_{t}\}$ and $\{\eta_{t}\}$ are independent, therefore 
expanding the term $A:= \cov\left\{h(Z_{t})K_{b}(Z_{t}-z)
\eta_{t},h(Z_{\tau})K_{b}(Z_{\tau}-z)\eta_{\tau}\right\}$ gives
\begin{displaymath}
A= \cov(\eta_{t},\eta_{\tau})\cdot 
\ex [h(Z_{t})K_{b}(Z_{t}-z)h(Z_{\tau})K_{b}(Z_{\tau}-z)].
\end{displaymath}
Since the covariance of the time series $\{\eta_{t}\}$ has size $\gu$, applying the Cauchy-Schwarz 
inequality gives
\begin{displaymath}
|A|\lesssim |t-\tau|^{-\gu}\cdot \ex |h(Z_{1})K_{b}(Z_{1}-z)|^2.
\end{displaymath}
Under Assumption   \ref{non:reg:ass:tec-part1}
the function $|h|^{2}\cdot f$ is uniformly bounded and $\|K\|_{2}<\infty$, therefore 
$\ex |h(Z_{1})K_{b}(Z_{1}-z)|^2\lesssim b^{-d}$, and hence
we obtain (\ref{a:p:non:reg:l:1:res3}). \hfill $\square$

\bigskip

\begin{lemma}\label{a:p:non:reg:l:2} 
Suppose the stationary
time series  $\{Z_{t}\}$ is $2$-mixing with size $\gv$ and Assumption \ref{non:reg:ass:tec-part1} is 
fulfilled for some $q,q_{F}\in(0,1)$.  If $1\leq t,\tau\leq T$, then
\begin{equation}\label{a:p:non:reg:l:2:res1}
|\cov\left\{K_{b}(Z_{t}-z),K_{b}(Z_{\tau}-z)\right\}|  
\lesssim\min\Bigl(b^{-d(1-q_{F})}; b^{-d(1+q)}|t-\tau|^{-q\mathfrak{v}}\Bigr).
\end{equation}
\end{lemma}

\bigskip

\noindent\textcolor{darkred}{\sc Proof.} The proof is very similar to the proof of Lemma \ref{a:p:non:reg:l:1} and we omit the details.\hfill $\square$

\bigskip

\begin{lemma}\label{a:p:non:reg:l:3} 
Suppose the assumptions in Theorem \ref{non:reg:theo-model1} are satisfied. Let 
$\widehat{g}$ be defined as in (\ref{non:reg:est:classical}). Then we have 
\begin{equation}\label{a:p:non:reg:l:3:res1}
\ex|\widehat{g}(z)-g(z)|\lesssim b^{2\rho}+b^{-d}T^{-1} +
b^{-d(1+q+q_{g}(1-[(q\gv)\vee1]))}T^{-[(q\gv)\wedge1]}+ b^{-d}T^{-(\gu\wedge1)}.
\end{equation}
\end{lemma}

\bigskip

\noindent\textcolor{darkred}{\sc Proof.} Consider the standard variance bias decomposition
\begin{equation}\label{a:p:non:reg:l:3:MSE:dec}
 \ex|\hat g(z)- g(z)|^2 = \var(\hat g(z)) + 
|\ex\hat g(z)-g(z)|^2.
\end{equation}
Under the stated assumptions we will derive the following two bounds, which altogether give
the estimate in (\ref{a:p:non:reg:l:3:res1}).  The bias is bounded by
\begin{gather}\label{a:p:non:reg:l:3:MSE:bias}
|\ex\hat g(z)-g(z)|^2\lesssim  b^{2\rho},
\end{gather}
while for the variance we have
\begin{gather}\label{a:p:non:reg:l:3:MSE:var}
\var(\hat g(z))\lesssim  b^{-d}T^{-1} +
b^{-d(1+q+q_{G}(1-[(q\gv)\vee1]))}T^{-[(q\gv)\wedge1]}+ b^{-d}T^{-(\gu\wedge1)}.
\end{gather}

We first prove (\ref{a:p:non:reg:l:3:MSE:bias}). We can write
\begin{displaymath}
\ex\hat g(z)= \frac{1}{T}\sum _{t=1}^{T}\ex\Bigl(\ex[X_{t}|Z_{t}]K_{b}(Z_{t}-z)\Bigr)=
\int\; du\; g(u) K_{b}(u-z).
\end{displaymath}
Since $g\in \mf{G}^{d}_{s,\triangle}$  and $K$ is a multiplicative kernel  
of order $r$ with $\int du |u|^r K(u)\leq S_{K}$, using a Taylor expansion 
up to the order $\rho=\min(r,s)$ leads to 
$\ex\hat g(z)=g(z)+b^{\rho} R$ with reminder 
$|R|\leq \triangle S_{K}<\infty,$ which proves  
(\ref{a:p:non:reg:l:3:MSE:bias}).

In order to proof 
(\ref{a:p:non:reg:l:3:MSE:var}), 
we consider the expansion
\begin{align}\nonumber
\var(\hat{g}(z)) &= \frac{1}{T^{2}} \sum_{t=1}^{T} 
\var\left\{ X_{t}K_{b}(Z_{t}-z)\right\}+ \frac{2}{T^{2}} 
\sum_{t>\tau} \cov\left\{X_{t}K_{b}(Z_{t}-z),X_{\tau}K_{b}(Z_{\tau}-z)\right\}\\
\label{a:p:non:reg:l:3:MSE:var:dec}
&=:A_{1}+A_{2}.
\end{align}
We will show that 
$ |A_{1}|\lesssim T^{-1} \cdot  b^{-d} + T^{-(\gu\wedge 1)} \cdot  b^{-d}$ and  
\begin{align}\label{a:p:non:reg:l:3:MSE:var:dec:A2}
|A_{2}|\lesssim\left\{
\begin{array}{ll}
  T^{-q\gv} \cdot b^{-d(1+q)}+ T^{-(\gu\wedge 1)} \cdot  b^{-d}, & q\mathfrak{v}\leq 1; \\
 T^{-1}\cdot   b^{-d} + T^{-1}\cdot b^{-d(1+q+q_{G}(1-q\mathfrak{v}))}+ T^{-(\gu\wedge 1)} \cdot  b^{-d}, & 
1<q\mathfrak{v}.
\end{array}
\right.
\end{align}
Furthermore, if $0 \leq q\mathfrak{v} \leq q/q_{G}+1$ then we show that
$|A_{1}|$ is dominated by $|A_{2}|$. 
Whereas for $q\mathfrak{v}>q/q_{G}+1$ 
the terms $|A_{1}|$ and $|A_{2}|$ are of the same order $O(T^{-1} \cdot  b^{-d} + T^{-(\gu\wedge 1)} 
\cdot  b^{-d})$. Therefore, the bounds derived for 
$|A_{2}|$ will lead to the estimates in 
(\ref{a:p:non:reg:l:3:MSE:var}). 

First let us consider $A_{1}$. Due to  stationarity of the process, 
we have the bound 
\begin{align*}
 T\cdot A_{1} \leq \ex|X_{1}K_{b}(Z_{1}-z)|^2&\lesssim   \ex|\varphi(Z_{1})K_{b}(Z_{1}-z)|^2 + 
\ex|h(Z_{1})K_{b}(Z_{1}-z)|^2.
\end{align*}  
Under the stated assumptions the functions $|\varphi|^p\cdot f$ with $p>2$ and $|h|^2\cdot f$ 
are uniformly  bounded and the kernel $\normV{K}_{2}<\infty$, therefore
$A_{1}\lesssim T^{-1}\cdot b^{-d}$. 

Let us now consider the  term $A_{2}$, which  is  bound by 
\begin{equation}
\label{a:p:non:reg:l:3:MSE:var:dec:cov}
T\cdot |A_{2}| \leq 
4\sum_{t=2}^T |\cov\left\{X_{t}K_{b}(Z_{t}-z),X_{1}K_{b}(Z_{1}-z)\right\}|,
\end{equation} 
where  using  representation (\ref{non:reg:model1}) the $t$-th summand in 
(\ref{a:p:non:reg:l:3:MSE:var:dec:cov}) can be estimated by
\begin{multline}
\label{a:p:non:reg:l:3:MSE:var:dec:cov:e1}
| \cov\left\{ \varphi(Z_{t})K_{b}(Z_{t}-z), \varphi(Z_{1})K_{b}(Z_{1}-z)\right\}|\\
\hfill+ |\cov\left\{ h(Z_{t})K_{b}(Z_{t}-z)\eta_{t}, h(Z_{1})K_{b}(Z_{1}-z)\eta_{1}\right\} |.
\end{multline}
If $q\mathfrak{v}\leq 1 $ and $\gu\leq 1$ then we 
bound the sum (\ref{a:p:non:reg:l:3:MSE:var:dec:cov})  using 
(\ref{a:p:non:reg:l:3:MSE:var:dec:cov:e1}), i.e., 
\begin{multline}
\label{a:p:non:reg:l:3:MSE:var:dec:cov:e2}
T\cdot |A_{2}|\lesssim \sum_{t=2}^T  | 
\cov\left\{ \varphi(Z_{t})K_{b}(Z_{t}-z), \varphi(Z_{1})K_{b}(Z_{1}-z)\right\}|\\
\hfill+\sum_{t=2}^T  |\cov\left\{ h(Z_{t})K_{b}(Z_{t}-z)\eta_{t}, h(Z_{1})K_{b}(Z_{1}-z)\eta_{1}\right\} |.
\end{multline}
We use the bounds (\ref{a:p:non:reg:l:1:res2}) and  (\ref{a:p:non:reg:l:1:res3}) 
in Lemma~\ref{a:p:non:reg:l:1} to estimate each of the sums in (\ref{a:p:non:reg:l:3:MSE:var:dec:cov:e2}) 
separately, which gives
\begin{equation}
\label{a:p:non:reg:l:3:MSE:var:dec:cov:e3}
T\cdot |A_{2}|\lesssim T^{-q\mathfrak{v} +1} \cdot b^{-d(1+q)}+T^{-\gu+1} \cdot b^{-d}.
\end{equation}
On the other hand if $q\mathfrak{v}> 1$ or if $\gu>1$ we partition the sum 
(\ref{a:p:non:reg:l:3:MSE:var:dec:cov}) into two parts, where  we estimate the first part using 
the bound (\ref{a:p:non:reg:l:1:res1})  in Lemma~\ref{a:p:non:reg:l:1} and the second using 
(\ref{a:p:non:reg:l:3:MSE:var:dec:cov:e1}), thus giving us 
\begin{multline}
\label{a:p:non:reg:l:3:MSE:var:dec:cov:e4}
T\cdot |A_{2}|  \lesssim h\cdot
 b^{-d(1-q_{G})}   + \sum_{t=h+1}^T  | \cov\left\{ \varphi(Z_{t})K_{b}(Z_{t}-z), 
\varphi(Z_{1})K_{b}(Z_{1}-z)\right\}|\\
\hfill+\sum_{t=h+1}^T  |\cov\left\{ h(Z_{t})K_{b}(Z_{t}-z)\eta_{t}, h(Z_{1})K_{b}(Z_{1}-z)\eta_{1}\right\} |.
\end{multline}
We use the bounds (\ref{a:p:non:reg:l:1:res2}) and  (\ref{a:p:non:reg:l:1:res3})  in 
Lemma~\ref{a:p:non:reg:l:1} to estimate each of the sums in (\ref{a:p:non:reg:l:3:MSE:var:dec:cov:e4}) 
separately, which gives
\begin{multline}\label{a:p:non:reg:l:3:MSE:var:dec:cov:e5}
T\cdot |A_{2}|  \lesssim h\cdot
 b^{-d(1-q_{G})}   +\left\{\begin{array}{ll} 
 T^{-q\gv+1} \cdot b^{-d(1+q)} + h^{-\gu+1} \cdot b^{-d},& q\gv\leq 1 \mbox{ and }\gu >1;\\
 h^{-q\gv+1} \cdot b^{-d(1+q)} + T^{-\gu+1} \cdot b^{-d},& q\gv> 1 \mbox{ and }\gu \leq1;\\
 h^{-q\gv+1} \cdot b^{-d(1+q)} + h^{-\gu+1} \cdot b^{-d},& q\gv> 1 \mbox{ and }\gu >1.\\
 \end{array}\right.
\end{multline}
Thereby using $h\approx b^{-d q_{G}}$ we obtain  
\begin{multline}
\label{a:p:non:reg:l:3:MSE:var:dec:cov:e6}
T\cdot |A_{2}|  \lesssim 
 b^{-d}   +\left\{\begin{array}{ll} 
 T^{-q\gv+1} \cdot b^{-d(1+q)} +  b^{-d(1+q_{G}(1-u))},& q\gv\leq 1 \mbox{ and }\gu >1;\\
 b^{-d(1+q+q_{G}(1-q\gv))} + T^{-\gu+1} \cdot b^{-d},& q\gv> 1 \mbox{ and }\gu \leq1;\\
b^{-d(1+q+q_{G}(1-q\gv))} +  b^{-d(1+q_{G}(1-u))},& q\gv> 1 \mbox{ and }\gu >1.\\
 \end{array}\right.
\end{multline}
and hence, combining (\ref{a:p:non:reg:l:3:MSE:var:dec:cov:e3}) and 
(\ref{a:p:non:reg:l:3:MSE:var:dec:cov:e6}) gives the  bound 
(\ref{a:p:non:reg:l:3:MSE:var:dec:A2}) for the term $A_{2}$.\hfill $\square$

\bigskip

We now state a slight variation of Theorem \ref{non:theo}, where $f$ can satisfy
slightly weaker conditions. We use this result to prove Theorem~\ref{non:reg:theo-model1}.
\begin{lemma}
\label{a:p:non:reg:l:4} 
Suppose the stationary time series $\{Z_{t}\}$ is $2$-mixing with size $\gv$ and 
Assumption \ref{non:reg:ass:tec-part1} is fulfilled for some $q,q_{F}\in(0,1)$.  
Let $\widehat{f}$ be defined as in (\ref{eq:density_estimator}), where  
the multiplicative kernel is of order $r>0$. In addition assume, that the function 
$f$ belongs to $\mf{G}^{d}_{s,\triangle}$ for $s,\triangle>0$ and let $\rho:=\min(r,s)$. Then we have
\begin{equation}
\label{a:p:non:reg:l:4:res1}
\ex|\widehat{f}(z)-f(z)|\lesssim 
b^{2\rho}+b^{-d}T^{-1}+b^{-d(1+q+q_{F}(1-[(q\gv)\vee1]))}T^{-[(q\gv)\wedge1]}.
\end{equation}
\end{lemma}

\bigskip

\noindent\textcolor{darkred}{\sc Proof.} Under the stated assumptions using 
Lemma \ref{a:p:non:reg:l:2} the proof  is very similar to the proof of Lemma \ref{a:p:non:reg:l:3} 
and we omit the details.\hfill $\square$

\bigskip

\noindent\textcolor{darkred}{\sc Proof of Theorem~\ref{non:reg:theo-model1}.} 
Consider the decomposition
\begin{align*}
\hat \varphi(z)-\varphi(z)&= 
\frac{\hat g(z)}{\hat f(z)}- 
\frac{\hat f(z)}{\hat f(z)}\varphi(z)\\
&= \frac{\hat g(z) - \hat f(z)\varphi(z)}{f(z)} +   
\frac {f(z)-\hat f(z)}{ \hat f(z)}\cdot   
\frac {\hat g(z)- \hat f(z)\varphi(z)}{  f(z)}.
\end{align*}
We first note that Lemma~\ref{a:p:non:reg:l:4} gives $\ex|f(z)-\hat f(z)|^2 =o(1)$,   
which implies that $|\hat f(z)^{-1}|$ is bounded in probability. Therefore the 
second term in the above expansion is of order $o_{P}(\{\hat g(z)- \hat f(z)\varphi(z)\}/f(z))$, 
hence in the decomposition above the  second term is negligible in comparison to the first term. 
Thereby bounding the first term of the decomposition we obtain the result. 
By using Lemma~\ref{a:p:non:reg:l:3} and \ref{a:p:non:reg:l:4} and noting that 
$q_{FG} = q_{F}\wedge q_{G}$, we obtain Theorem~\ref{non:reg:theo-model1}. 
\hfill $\square$

\bigskip

\noindent\textcolor{darkred}{\sc Proof of Corollary~\ref{non:reg:cor-model1}}
Under the assumption on the bandwidth the result is obtained by balancing the terms 
in the bound given in Theorem~\ref{non:reg:theo-model1}.\hfill $\square$

\bigskip

\begin{lemma}\label{a:p:non:reg:l:5} 
Suppose the stationary
vector time series $\{(X_{t},Z_{t})\}$ is $2$-mixing with size $\gv$ and Assumption 
\ref{non:reg:ass:tec-part2} is fulfilled for some $q,q_{G}\in(0,1)$.  If $1\leq t,\tau\leq T$, then
\begin{equation}\label{a:p:non:reg:l:5:res1}
|\cov\left\{X_{t}K_{b}(Z_{t}-z),X_{\tau}K_{b}(Z_{\tau}-z)\right\}|
\lesssim \min \Bigl(b^{-d(1-q_{G})}, b^{-d(1+q)} |t-\tau|^{-q\gv}\Bigr).\end{equation}
\end{lemma}

\bigskip

\noindent\textcolor{darkred}{\sc Proof.} Under Assumption \ref{non:reg:ass:tec-part2} (ii) the 
first bound in (\ref{a:p:non:reg:l:5:res1}) 
follows from (\ref{a:p:non:reg:l:1:res1}) in  Lemma~\ref{a:p:non:reg:l:1:res1}. On the other hand, 
under Assumption   \ref{non:reg:ass:tec-part1} (i,iii) the function 
$\ex[|X_{1}|^{p}|Z_{1}]\cdot f$ is uniformly bounded and  $\|K\|_{p}$ is finite  for some 
$p=2/(1-q)>2$, therefore we have 
$[\ex |X_{1}K_{b}(Z_{1}-z)|^p]^{2/p}\lesssim b^{-d(q+1)}.$ Using the $2$-mixing property of 
$\{Z_{t}\}$ together with  \cite{b:hal-hey-80}, Theorem A.6, we obtain  the second bound in  
(\ref{a:p:non:reg:l:5:res1}). \hfill $\square$

\bigskip

\begin{lemma}\label{a:p:non:reg:l:6} 
Suppose the stationary 
vector time series $\{(X_{t},Z_{t})\}$ is $2$-mixing with size $\gv$ and Assumption
\ref{non:reg:ass:tec-part2} is  fulfilled for some $q,q_{G}\in(0,1)$.  
Let $\widehat{g}$ be defined as in (\ref{non:reg:est:classical}), where the 
multivariate kernel is of order $r>0$. In addition assume, that the function $g=\varphi\cdot f$ 
belongs to $\mf{G}^{d}_{s,\triangle}$ for $s,\triangle>0$ and let $\rho:=\min(r,s)$. Then we have
\begin{equation}
\label{a:p:non:reg:l:6:res1}
\ex|\widehat{g}(z)-g(z)|\lesssim b^{2\rho}+b^{-d}T^{-1}+b^{-d(1+q+q_{G}(1-[(q\gv)\vee1]))}
T^{-[(q\gv)\wedge1]}
\end{equation}
\end{lemma}

\bigskip

\noindent\textcolor{darkred}{\sc Proof.}  Under the stated assumptions using 
Lemma \ref{a:p:non:reg:l:5} the proof is very similar to the proof of Lemma \ref{a:p:non:reg:l:3} 
and we omit the details. \hfill $\square$

\bigskip

\noindent\textcolor{darkred}{\sc Proof of Theorem~\ref{non:reg:theo-model2}.} Using  
Lemma~\ref{a:p:non:reg:l:4} and \ref{a:p:non:reg:l:4} we obtain the result using a similar 
proof as Theorem~\ref{non:reg:theo-model1}.\hfill $\square$

\bigskip

\noindent\textcolor{darkred}{\sc Proof of Corollary~\ref{non:reg:cor-model2}}
Under the assumption on the bandwidth the result is obtained by balancing the terms in the 
bound given in Theorem~\ref{non:reg:theo-model1}.\hfill $\square$

%% file: appendix-proofs-nonpar-panel.tex
\subsection{Proofs: Nonparametric regression 
 for panel time series}\label{a:p:non:pan}

\begin{lemma}\label{a:p:non:pan:l:cov:g} 
Suppose $\{X_{t,i}\}$ satisfies (\ref{non:pan:model}),  
for all $i,j\in \N$, $\{(X_{t,i},Z_{t,i},X_{t,j},Z_{t,j})\}$ is a stationary time series, 
and the panel time series $\{(X_{t,i},Z_{t,i})\}$ is $2$-mixing with 
size $\mathfrak{v}$ and $\mathfrak{u}$ (as defined in Definition \ref{non:pan:def:t-dep}) and 
Assumption \ref{non:pan:ass:tec} is satisfied for some $ q_{G},q\in(0,1)$. 
If $1\leq t,\tau\leq T$ and $1\leq i<j \leq N$, then
\begin{multline}\label{a:p:non:pan:l:cov:g:ij}
 |\cov\left\{X_{t,i}K_{b_{i}}(Z_{t,i}-z),X_{\tau,j}K_{b_{j}}(Z_{\tau,j}-z)
\right\}|\lesssim  \\
\hfill  \min\Bigl((b_{i}b_{j})^{-\frac{d}{2}(1-q_{G})}; (b_{i}b_{j})^{-\frac{d}{2}(q+1)} 
|t-\tau|^{-q\mathfrak{u}} \Bigr),
\end{multline}
while if $1\leq t,\tau\leq T$ and $1\leq i \leq N$, then
\begin{multline}\label{a:p:non:pan:l:cov:g:j}
|\cov\left\{X_{t,i}K_{b_{i}}(Z_{t,i}-z),X_{\tau,i}K_{b_{i}}(Z_{\tau,i}-z)\right\}|\lesssim
\\
\hfill \min\Bigl(b_{i}^{-d(1-q_{G})}; 
b_{i}^{-d (q+1)} |t-\tau|^{-q\mathfrak{v}} \Bigr).
\end{multline}
\end{lemma}

\bigskip

\noindent\textcolor{darkred}{\sc Proof.} 
Using Assumption \ref{non:pan:ass:tec} together with H\"older's inequality, and recalling that
$q_{G} = 1 - 2/p_{G}$  with $p_{G}^{-1}+\bar p_{G}^{-1}=1$ and  
$\|G_{t,\tau}^{(i,j)}\|_{p_{G}}$ is uniformly bounded,  we have 
\begin{displaymath}
|\cov\left\{ X_{t,i}K_{b_{i}}(Z_{t,i}-z),
X_{\tau,j}K_{b_{j}}(Z_{\tau,j}-z)\right\} | 
\lesssim
(b_{i}\cdot b_{j})^{-d/p_{G}},
\end{displaymath}
where the bound is obtained by using  Lyaponov's inequality, which gives $\|K\|_{p_{G}}<\infty$. 
This gives the common  bound 
in  (\ref{a:p:non:pan:l:cov:g:ij}) and (\ref{a:p:non:pan:l:cov:g:j}). On the other hand,  
we have $\ex[|X_{t,i}K_{b_{i}}(Z_{t,i}-z)|^{p}
]<\infty$ with $p=2/(1-q)>2$ (Assumption \ref{non:pan:ass:tec} (i)). Therefore, 
using the $2$-mixing property of the panel time series
$\{(X_{t,i},Z_{t,i})\}$  together with  \cite{b:hal-hey-80}, Theorem A.6, for 
$i\ne j$, we obtain
\begin{multline}\label{a:p:non:pan:l:cov:g:e1}
|\cov\left\{ X_{t,i}K_{b_i}(Z_{t,i}-z),
X_{\tau,j}K_{b_{j}}(Z_{\tau,j}-z)\right\} | \\
\hfill
\lesssim\{\ex[|X_{1,i}K_{b_{i}}(Z_{1,i}-z)|^{p} ]
\cdot\ex[|X_{1,j}K_{b_{j}}(Z_{1,j}-z)|^{p} ]\}^{1/p} \cdot |t-\tau|^{-q\mathfrak{u}},
\end{multline}
while for $i=j$
\begin{multline}\label{a:p:non:pan:l:cov:g:e2}
|\cov\left\{ X_{t,i}K_{b_i}(Z_{t,i}-z),
X_{\tau,i}K_{b_{i}}(Z_{\tau,i}-z)\right\} | \\
\hfill
\lesssim\{\ex[|X_{1,i}K_{b_{i}}(Z_{1,i}-z)|^{p} ]\}^{2/p} \cdot |t-\tau|^{-q\mathfrak{v}}.
\end{multline}
Since under Assumption~\ref{non:pan:ass:tec}, the function
$g_{i}^{(p)}(\cdot)=\ex[|X_{1,i}|^{p}|Z_{1,i}=\cdot]f_{i}(\cdot)$ is uniformly bounded 
and $\|K\|_{p}<\infty$ we have
\begin{align}\label{a:p:non:pan:l:cov:g:e3}
  &\ex[|X_{1,i}K_{b_{i}}(Z_{1,i}-z)|^{p} ]^{1/p} \lesssim b_{i}^{-\frac{d}{2} (q+1)}.
\end{align} 
Therefore, (\ref{a:p:non:pan:l:cov:g:e1}) together with
(\ref{a:p:non:pan:l:cov:g:e3}) gives the second bound in
(\ref{a:p:non:pan:l:cov:g:ij}), where (\ref{a:p:non:pan:l:cov:g:e2}) and
(\ref{a:p:non:pan:l:cov:g:e3}) leads to the second bound in
(\ref{a:p:non:pan:l:cov:g:j}), which proves the result.\hfill $\square$

\bigskip

\begin{lemma}\label{a:p:non:pan:l:cov:f}  
Suppose $\{X_{t,i}\}$ satisfies (\ref{non:pan:model}),  
for all $i,j\in \N$, $\{(X_{t,i},Z_{t,i},X_{t,j},Z_{t,j})\}$ is a stationary time series, 
and the panel time series $\{(X_{t,i},Z_{t,i})\}$ is $2$-mixing with 
size $\mathfrak{v}$ and $\mathfrak{u}$ (as defined in Definition \ref{non:pan:def:t-dep}) and 
Assumption \ref{non:pan:ass:tec} is satisfied for some $q_{F},q\in(0,1)$. 
If $1\leq t,\tau\leq T$ and $1\leq i<j \leq N$, then
\begin{equation}
\label{a:p:non:pan:l:cov:f:ij}
 |\cov\left\{K_{b_{i}}(Z_{t,i}-z),K_{b_{j}}(Z_{\tau,j}-z)\right\}|\lesssim 
\min\Bigl((b_{i}b_{j})^{-\frac{d}{2}(1-q_{F})}; (b_{i}b_{j})^{-\frac{d}{2}(1+q)} |t-\tau|^{-\mathfrak{u}} \Bigr),
\end{equation}
while if $1\leq t,\tau\leq T$ and $1\leq i \leq N$, then
\begin{equation}\label{a:p:non:pan:l:cov:f:j}
|\cov\left\{K_{b_{i}}(Z_{t,i}-z),K_{b_{i}}(Z_{\tau,i}-z)\right\}|  \lesssim
\min\Bigl(b_{i}^{-d(1-q_{F})}; 
b_{i}^{-d(1+q) }|t-\tau|^{-\mathfrak{v}}
\Bigr).
\end{equation}
\end{lemma}

\bigskip

\noindent\textcolor{darkred}{\sc Proof.}  The proof is very similar to the proof of 
Lemma \ref{a:p:non:pan:l:cov:g} and we omit the details.\hfill $\square$

We use the lemma below to prove Theorem \ref{non:pan:theo:MSE}, which requires the 
following definitions
\begin{eqnarray*}
g := \frac{1}{N}\sum_{i=1}^N g_{i}, \quad  
\hat g := \frac{1}{N}\sum_{i=1}^N \hat g_{i}, \quad f := \frac{1}{N}\sum_{i=1}^N f_{i} \quad \textrm{and}\quad
\hat f := \frac{1}{N}\sum_{i=1}^N \hat f_{i}.
\end{eqnarray*}

\begin{lemma}\label{a:p:non:pan:l:MSE} 
Let us suppose that all assumptions in Theorem \ref{non:pan:theo:MSE} hold. Then we have 
\begin{multline}\label{a:p:non:pan:l:MSE:g}
 \ex|\hat g(z)- g(z)|^2\lesssim  \frac{1}{N}
\sum_{i=1}^N\Bigl\{ b_{i}^{2\rho_{i}} +T^{-1}\cdot  b_{i}^{-d}+   T^{-[(q\gu)\wedge 1]} \cdot 
 b_{i}^{-d(1+q+q_{G}-q_{G}[(q\mathfrak{u})\vee 1])}\\
\hfill+
N^{-1}\cdot T^{-[(q\gv)\wedge 1]} \cdot 
 b_{i}^{-d(1+q+q_{G}-q_{G}[(q\mathfrak{v})\vee 1])}\Bigr\}\end{multline}
and
\begin{multline}\label{a:p:non:pan:l:MSE:f}
 \ex|\hat f(z)- f(z)|^2\lesssim  \frac{1}{N}
\sum_{i=1}^N\Bigl\{ b_{i}^{2\rho_{i}} +T^{-1}\cdot  b_{i}^{-d}+T^{-[(q\gu)\wedge 1]} \cdot 
 b_{i}^{-d(1+q+q_{F}-q_{F}[(q\mathfrak{u})\vee 1])}\\
\hfill+
N^{-1}\cdot T^{-[(q\gv)\wedge 1]} \cdot 
b_{i}^{-d(1+q+q_{F}-q_{F}[(q\mathfrak{v})\vee 1])}\Bigr\}.\end{multline}
\end{lemma}
\bigskip

\noindent\textcolor{darkred}{\sc Proof.} We only give the details for the proof of the MSE of 
the estimator $\hat g$. The proof of the other result is very similar (but uses the bounds given in Lemma \ref{a:p:non:pan:l:cov:f} rather than Lemma \ref{a:p:non:pan:l:cov:g}) and we omit the 
details. Consider the standard variance and  bias decomposition
\begin{equation}
 \ex|\hat g(z)- g(z)|^2 = \var(\hat g(z)) + 
|\ex\hat g(z)-g(z)|^2.
\end{equation}
Under the stated assumptions we will derive the following two bounds. 
The bias is bounded by
\begin{gather}\label{a:p:non:pan:l:MSE:g:eq:bias}
|\ex\hat g(z)-g(z)|^2\lesssim  \frac{1}{N}
\sum_{i=1}^N b_{i}^{2\rho_{i}},
\end{gather}
while for the variance we have 
\begin{multline}\label{a:p:non:pan:l:MSE:g:eq:var}
\var(\hat g(z))\lesssim T^{-[(q\gu)\wedge 1]} \cdot 
\frac{1}{N}\sum_{i=1}^N[   b_{i}^{-d} + b_{i}^{-d(1+q+q_{G}-q_{G}[(q\mathfrak{u})\vee 1])}]\\
\hfill+
(NT^{[(q\gv)\wedge 1]})^{-1} \cdot 
\frac{1}{N}\sum_{i=1}^N[   b_{i}^{-d} + b_{i}^{-d(1+q+q_{G}-q_{G}[(q\mathfrak{v})\vee 1])}].
\end{multline}
We now prove (\ref{a:p:non:pan:l:MSE:g:eq:bias}). Using iterative conditional  
expectation we can write
\begin{displaymath}
\ex\hat g(z)= \frac{1}{NT}\sum _{i=1}^{N}
\sum _{t=1}^{T}\ex\Bigl(\ex[X_{t,i}|Z_{t,i}]K_{b_{i}}(Z_{t,i}-z)\Bigr)=
\frac{1}{N}\sum _{i=1}^{N}\int\; du\; g_{i}(u) K_{b_{i}}(u-z)
\end{displaymath}
where $g_{i}= \ex[X_{t,i}|Z_{t,i}=\cdot]f_{i}(\cdot)$ for all $t,i$. 
Since $g=\frac{1}{N}\sum_{i=1}^N g_{i}$ with 
$g_{i}\in \mf{G}^{d}_{s_{i},\triangle}$  and $K$ is a multiplicative kernel  
of order $r$ with $\int du |u|^r K(u)\leq S_{K}$, using a Taylor expansion 
up to order  $\rho_{i}=\min(r,s_{i})$ leads to 
$\ex\hat g(z)=g(z)+\frac{1}{N}\sum _{i=1}^{N}b_{i}^{\rho_{i}} R_{i}$ with reminder 
$|R_{i}|\leq \triangle S_{K}<\infty.$ Thus applying Jensens inequality we obtain 
(\ref{a:p:non:pan:l:MSE:g:eq:bias}).

In order to proof 
(\ref{a:p:non:pan:l:MSE:g:eq:var}) we consider the expansion
\begin{align}\label{a:p:non:pan:l:MSE:g:eq:var:dec}
&\var(\hat{g}(z)) = A_{1}+A_{2}+A_{3}+A_{4}\intertext{with}\nonumber
&A_{1}=\frac{1}{N^2T^{2}} \sum_{t=1}^{T}\sum_{i=1}^{N} 
\var\left\{ X_{t,i}K_{b_{i}}(Z_{t,i}-z)\right\}\;,\\\nonumber
&A_{2}= \frac{2}{N^2T^{2}} \sum_{t=1}^{T}\sum_{j>i}
\cov\left\{X_{t,i}K_{b_{i}}(Z_{t,i}-z),
X_{t,j}K_{b_{j}}(Z_{t,j}-z)\right\}\;, \\\nonumber
&A_{3} = \frac{4}{N^2T^{2}} \sum_{t>\tau}\sum_{j>i} 
\cov\left\{X_{t,i}K_{b_{i}}(Z_{t,i}-z),X_{\tau,j}K_{b_{j}}(Z_{\tau,j}-z)\right\}\;, 
\\\nonumber
&A_{4}  =\frac{2}{N^{2}T^{2}} \sum_{t>\tau}\sum_{i=1}^{N} 
\cov\left\{X_{t,i}K_{b_{i}}(Z_{t,i}-z),X_{\tau,i}K_{b_{i}}(Z_{\tau,i}-z)\right\}.
\end{align}
We will show that $|A_{1}|,|A_{2}| \lesssim T^{-1} \cdot  \frac{1}{N}
\sum_{i=1}^N  b_{i}^{-d}$, 
\begin{align*}
 |A_{3}|&\lesssim T^{-[(q\gu)\wedge 1]} \cdot 
\frac{1}{N}\sum_{i=1}^N[   b_{i}^{-d} + b_{i}^{-d(1+q+q_{G}-q_{G}[(q\mathfrak{u})\vee 1])}]\mbox{ and }\\
 |A_{4}|&\lesssim 
(NT^{[(q\gv)\wedge 1]})^{-1} \cdot 
\frac{1}{N}\sum_{i=1}^N[   b_{i}^{-d} + b_{i}^{-d(1+q+q_{G}-q_{G}[(q\mathfrak{v})\vee 1])}]. \end{align*}
Furthermore, if $0 \leq q(\mathfrak{v}\wedge \gu) \leq q/q_{G}+1$ then the terms $|A_{1}|$ and $|A_{2}|$ are
dominated by $|A_{3}|+|A_{4}|$. Whereas for $q(\mathfrak{v}\wedge \gu)>q/q_{G}+1$ all
the terms are of the same order. Therefore, the bound derived for 
$|A_{3}|+|A_{4}|$ will lead to the estimate in (\ref{a:p:non:pan:l:MSE:g:eq:var}). 

First let us consider $A_{1}$. Due to stationarity, 
we obtain the bound 
\begin{align*}
N\cdot T\cdot A_{1} \leq 
\frac{1}{N}\sum_{i=1}^N\ex[X_{1,i}^2K_{b_{i}}^2(Z_{1,i}-z)].
\end{align*}  
Thereby, under the assumption that the functions $g_i^{(2)}(\cdot):= 
\ex[|X_{1,i}|^2|Z_{1,i}=\cdot]f_{i}(\cdot)$ are uniformly bounded and the kernel $\|K\|_{2}<\infty$,  
we have $A_{1}\lesssim (N\cdot T)^{-1} \cdot \frac{1}{N}\sum_{i=1}^N  b_{i}^{-d}$. 

It is straightforward to show that $T\cdot |A_{2}|$ 
is bounded by
\begin{displaymath}
\frac{2}{N^2}\sum_{j>i} \var(X_{1,i}K_{b_{i}}(Z_{1,i}-z))^{1/2}\var(X_{1,j}K_{b_{j}}(Z_{1,j}-z))^{1/2} 
\lesssim \frac{1}{N^2}\sum_{j>i} (b_{i}b_{j})^{-d/2},\end{displaymath}
where the inequality  above follows by applying the same 
arguments as those used for $A_{1}$. Therefore using  the Cauchy-Schwarz inequality 
we obtain  $|A_{2}|\lesssim  T^{-1}\cdot 
\frac{1}{N}\sum_{i=1}^N b_{i}^{-d}$.

The term $T\cdot |A_{3}|$ is bounded by the sum
\begin{displaymath} 
\frac{8}{N^2}\sum_{j>i}\sum_{t=2}^T |\cov\left\{X_{t,i}K_{b_{i}}(Z_{t,i}-z),X_{1,j}
K_{b_{j}}(Z_{1,j}-z)\right\}|.
\end{displaymath} We now derive bounds for
$ T\cdot |A_{3}|$ for different mixing sizes. If $q\gu \leq 1$ then we estimate the inner sum using the second bound in 
(\ref{a:p:non:pan:l:cov:g:ij}) of Lemma~\ref{a:p:non:pan:l:cov:g}, i.e., $T\cdot |A_{3}|  \lesssim \frac{1}{N^2}\sum_{j>i} (b_{i}b_{j})^{-\frac{d}{2}(q+1)} T^{-q\gu+1}$, which leads together with 
the Cauchy-Schwarz inequality to $|A_{3}|  \lesssim T^{-q\gu}\frac{1}{N}\sum_{i=1}^N b_{i}^{-d(q+1)} $. On the other hand, if $q\gu > 1$ then we partition the inner  sum into two parts which we estimate separately using now the two bounds in 
(\ref{a:p:non:pan:l:cov:g:ij}) of Lemma~\ref{a:p:non:pan:l:cov:g}, thus giving us
\begin{align*}
T\cdot |A_{3}|  \lesssim&\frac{1}{N^2}\sum_{j>i} \Bigl\{\sum_{t=2}^{h_{ij}}
 (b_{i}b_{j})^{-\frac{d}{2}(1-q_{G})} 
+
\sum_{t=h_{ij}+1}^T (b_{i}b_{j})^{-\frac{d}{2}(q+1)} 
t^{-q\mathfrak{u}}\Bigr\}.
\end{align*}
Thereby
using $h_{ij}\approx (b_{i}b_{j})^{-\frac{d}{2}q_{G}}$ together with the Cauchy-Schwarz inequality we obtain 
$T\cdot |A_{3}|\lesssim\frac{1}{N}\sum_{i=1}^N[b_{i}^{-d}+ b_{i}^{-d(1+q+q_{G}(1-q\mathfrak{u}))} ].$ Combining the bounds in the two cases $q\gu \leq 1$ and $q\gu > 1$ we have  
$|A_{3}|\lesssim T^{-[(q\gu)\wedge 1]}\cdot \frac{1}{N}\sum_{i=1}^N[b_{i}^{-d}+ b_{i}^{-d(1+q+q_{G}-q_{G}[(q\mathfrak{u})\vee1])}].$

The term $N\cdot T\cdot |A_{4}|$ is  bounded  by 
\begin{displaymath}\frac{4}{N}\sum_{i=1}^N\sum_{t=2}^T |\cov\left\{X_{t,i}K_{b_{i}}(Z_{t,i}-z),
X_{1,i}K_{b_{i}}(Z_{1,i}-z)\right\}|.\end{displaymath} 
We estimate the inner sum applying the same arguments as those used for $|A_{3}|$ (but use the bounds given  in (\ref{a:p:non:pan:l:cov:g:j})  rather than (\ref{a:p:non:pan:l:cov:g:ij}) of 
Lemma~\ref{a:p:non:pan:l:cov:g}). Thereby we obtain $N\cdot |A_{4}|\lesssim T^{-[(q\gv)\wedge 1]}\cdot \frac{1}{N}\sum_{i=1}^N[b_{i}^{-d}+ b_{i}^{-d(1+q+q_{G}-q_{G}[(q\mathfrak{v})\vee1])}].$\hfill $\square$

\bigskip

\noindent\textcolor{darkred}{\sc Proof of Theorem~\ref{non:pan:theo:MSE}.} The proof is very similar to the 
proof of Theorem~\ref{non:reg:theo-model1} (but uses Lemma~\ref{a:p:non:pan:l:MSE}  rather than Lemma~\ref{a:p:non:reg:l:3} and \ref{a:p:non:reg:l:4}) and we omit the details.\hfill $\square$

\bigskip

\noindent\textcolor{darkred}{\sc Proof of Corollary~\ref{non:pan:cor:MSE:gu}.}
Under the assumptions of the corollary we have
\begin{equation}\label{a:p:non:pan:cor:MSE:gu:e} |\hat \varphi(z)- \varphi(z)|^{2} = 
O_{p} \Bigl(\frac{1}{N}
\sum_{i=1}^N\Bigl\{ b_{i}^{2\rho_{i}} +T^{-1}\cdot b^{-d}+ T^{-[(q\gu)\wedge 1]} \cdot 
 b_{i}^{-d(1+q+q_{FG}-q_{FG}[(q\mathfrak{u})\vee 1])}\Bigr\}\Bigr),
\end{equation}
applying Theorem~\ref{non:pan:theo:MSE},  where  $\frac{1}{N}
\sum_{i=1}^N N^{-1}\cdot T^{-[(q\gv)\wedge 1]} \cdot 
 b_{i}^{-d(1+q+q_{FG}-q_{FG}[(q\mathfrak{v})\vee 1])}$ is asymptotically negligible when $\gv\geq \gu$. Under the assumption on the bandwidths the result is obtained   balancing  the three terms in the bound  (\ref{a:p:non:pan:cor:MSE:gu:e}).\hfill $\square$

\bigskip

\noindent\textcolor{darkred}{\sc Proof of Corollary~\ref{non:pan:cor:MSE:gr}.}
The result is obtained applying Theorem~\ref{non:pan:theo:MSE}, where given $\gv\leq \gu$ the  assumption on the bandwidths provides the balance of the four terms of the bound (\ref{non:pan:theo:MSE:e}).\hfill $\square$

\bigskip

\noindent\textcolor{darkred}{\sc Proof of Corollary~\ref{non:pan:cor:MSE:gr:N}.}
Under the assumptions of the corollary we apply  Corollary~\ref{non:pan:cor:MSE:gr}, where given $N^{\zeta_{i}}\approx T^{(\gamma_{i}-\delta_{i})}$ for some $\zeta_{i}> 0$ the  bound (\ref{non:pan:cor:MSE:gr:e}) simplifies to 
\begin{equation*}\label{a:p:non:pan:cor:MSE:gr:N:e}
  |\hat \varphi(z)- \varphi(z)|^{2}=O_{P}\Bigl(\frac{1}{N}
\sum_{i=1}^N [T^{\{ \frac{\gamma_{i}-\delta_{i}}{\gamma_{i}}-\frac{\gamma_{i}-\delta_{i}}{\delta_{i}}\}\cdot[(q\gu)\wedge1]}+T^{-\frac{\gamma_{i}-\delta_{i}}{\zeta_{i}}+\frac{\gamma_{i}-\delta_{i}}{\delta_{i}}\cdot[(q\gv)\wedge1]}]\cdot  T^{-\frac{2\rho_{i}}{2\rho_{i}+d}\cdot\gamma_{i}}\Bigr). \end{equation*}
Since $\gv\leq \gu$ for each $i\in \N$  implies $T^{\{ \frac{\gamma_{i}-\delta_{i}}{\gamma_{i}}-\frac{\gamma_{i}-\delta_{i}}{\delta_{i}}\}\cdot[(q\gu)\wedge1]}=O(1)$. We obtain the result, if  for each $i\in \N$ we have   $T^{-\frac{\gamma_{i}-\delta_{i}}{\zeta_{i}}+\frac{\gamma_{i}-\delta_{i}}{\delta_{i}}\cdot[(q\gv)\wedge1]}=O(1)$ or equivalently  $\delta_{i}/[(q\gv)\wedge 1]\geq \zeta_{i}$, which is just the condition given in Corollary~\ref{non:pan:cor:MSE:gr:N}.\hfill $\square$

%% file: appendix-proofs-cov.tex
\subsection{Covariances and 2-mixing rates for linear
processes}\label{sec:cov-2mix}
We use the results derived in this section in Section \ref{sec:linear}, 
where we compared the rates of convergence for linear processes
with the rates in the general 2-mixing case. 

Let us suppose $\{Z_{t}\}$ satisfies the linear process representation 
in (\ref{eq:MAinfty}). By placing some additional 
conditions on the innovations we have the 
following lemma, which is due to \cite{p:gir-kou-96}, Lemma 1 and 2.
\begin{lemma}[\cite{p:gir-kou-96}]
\label{lemma:MAdensity}
Suppose $\{Z_{t}\}$ is a linear process which satisfies (\ref{eq:MAinfty}), and 
$\cov(Z_{0},Z_{t})\leq Ct^{-\theta}$
Let $f$ be the density of $Z_{t}$ and $f_{t}$ denote the joint density
$Z_{0},Z_{t}$.
If $\Ex(|\varepsilon_{t}^{3}|)<\infty$, and for all $u\in \mathbb{R}$ suppose the 
characteristic function satisfies $|\Ex[\exp(-iu\varepsilon_{1})]|\leq \frac{1}{(1+|u|)^{\delta}}$ 
for some $\delta > 0$,  then the joint density satisfies the relation
\begin{eqnarray*}
f_{t}(x,y) = f(x)f(y) + r(t)f^{\prime}(x)f^{\prime}(y) + O(t^{-\theta-d}),
\end{eqnarray*} 
where $f^{\prime}\in L_{1}(\mathbb{R})$ and $r(t)=\cov(Z_{0},Z_{t})$, for 
some $0< d < \min(\frac{\theta}{7},\frac{1-\theta}{12})$.
\end{lemma}
Using the result above the MSE of the kernel estimator with observations
from a linear process can be derived. 

For most processes, there isn't a direct correspondence between 
the 2-mixing  and the covariance size. However
for Gaussian processes both sizes are linked by the inequality
\begin{equation}\label{eq:cov-mix}
\frac{|\cov(X_{0},X_{t})|}{\var(X_{0})} \leq 
\sup_{A\in \sigma(Z_{0}),B\in
\sigma(Z_{t})} |P(A\cap B) - P(A)P(B)| \leq 2\pi \frac{|\cov(X_{0},X_{t})|}{\var(X_{0})}
\end{equation}
(see \cite{b:dou-94}, Section 2.1), thus the 
covariance and the 2-mixing sizes are the same. 
Suppose that $\{Z_{t}\}$ satisfies (\ref{eq:MAinfty}), where the innovations
are Gaussian and $|a_{j}|\lesssim j^{-\theta}$. Then we have 
\begin{equation}\label{eq:cov}
\left.
\begin{array}{c}
\frac{|\cov(X_{0},X_{t})|}{\var(X_{0})} \\
\textrm{ and } \\
\sup\limits_{A\in \sigma(Z_{0})
B\in\sigma(Z_{t})} |P(A\cap B) - P(A)P(B)|
\end{array}
\right\}
\lesssim 
\left\{
\begin{array}{cc}
(\log t) t^{-(2\theta-1) }, & \textrm{ if }  1/2 < \theta \leq  1; \\
t^{-\theta},              & \textrm{ if } \theta > 1.
\end{array}
\right.
\end{equation}

We now consider more general linear processes, which are not
necessarily Gaussian. Then the   
covariance size does not immediately give the 2-mixing size. 
However, if the density of the innovations satisfies certain smoothness conditions
then we can obtain the following bound. 
\begin{lemma}\label{lemma:MA}
Suppose $\{Z_{t}\}$ is a linear process which satisfies the representation
$Z_{t} = \sum_{j=0}^{\infty} a_{j}\varepsilon_{t-j}$, where the parameters
$|a_{j}|\leq C j^{-\theta}$ and $\theta > 1/2$. 
Let $f_{\varepsilon}$ be the density of the innovation $\varepsilon_{t}$. 
If $\Ex(|\varepsilon_{t}|^{\ell})<\infty$ (where $\ell>2$) and $\int
|f_{\varepsilon}(x+a)-f_{\varepsilon}(x)|dx\leq C|a|$,  then we have
\begin{eqnarray*}
\sup_{A\in \sigma(Z_{0}),B\in\sigma(Z_{t})} |P(A\cap B) - P(A)P(B)| \leq 
C j^{(-2\theta + 1)\frac{\ell}{2(\ell+1)}}
\end{eqnarray*}
where $C$ are some arbitrary constants.
\end{lemma}
PROOF. The result can be proved using a straightforward adaptation of the proofs in \cite{p:cha-74}, 
\cite{p:gor-77} and \cite{b:dav-94} (Theorem 14.9), who proved the result for 
strong $\alpha$-mixing. Hence we omit the details. \hfill $\Box$

\begin{remark}
It is interesting to compare the 2-mixing sizes derived in Lemma \ref{lemma:MA} with
the strong $\alpha$-mixing results 
for MA$(\infty)$ processes. Under the same set of conditions, but with the additional
restriction that $\theta > 3/2$, we have that 
$$\sup_{A\in \sigma(Z_{0},Z_{-1},\ldots),B\in\sigma(Z_{t},Z_{t+1},\ldots)} 
|P(A\cap B) - P(A)P(B)| \lesssim |t|^{(-2\theta + 1)\frac{\ell}{2(\ell+1)+1}}.$$ 
In other words, the 2-mixing size is larger  than the $\alpha$-mixing size. This is because, by definition,
the $\sigma$-algebras involved in the definition of  $\alpha$-mixing is far
larger than the $\sigma$-algebras in the definition of 2-mixing, thus 
allowing more extreme cases.\hfill $\Box$
\end{remark}
Comparing Lemma \ref{lemma:MA} with the covariance size given in (\ref{eq:cov})
we see when the Gaussianity assumption is relaxed the covariance  and 
2-mixing sizes no longer coincide. However by using Lemma \ref{lemma:MA} and 
\cite{b:hal-hey-80}, Theorem A.5, we have the upper and lower bounds
\begin{eqnarray*}
j^{(-2\theta + 1)\frac{\ell}{(\ell-2)}} 
\lesssim  \sup_{A\in \sigma(Z_{0}),B\in\sigma(Z_{t})} |P(A\cap B) - P(A)P(B)| \lesssim 
j^{(-2\theta + 1)\frac{\ell}{2(\ell+1)}}. 
\end{eqnarray*}
Therefore the 2-mixing size $\gv$ of the linear process $\{Z_{t}\}$  is bounded by 
\begin{eqnarray}\label{eq:up:lo:bd}
(2\theta - 1)\frac{\ell}{2(\ell+1)}\leq \gv \leq (2\theta - 1)\frac{\ell}{(\ell-2)}.
\end{eqnarray}

%% file: nonparametric.bbl
\begin{thebibliography}{39}
\providecommand{\natexlab}[1]{#1}
\providecommand{\url}[1]{\texttt{#1}}
\expandafter\ifx\csname urlstyle\endcsname\relax
  \providecommand{\doi}[1]{doi: #1}\else
  \providecommand{\doi}{doi: \begingroup \urlstyle{rm}\Url}\fi

\bibitem[Athreya and Pantula(1986)]{p:ath-pan-86}
K.~B. Athreya and S.~Pantula.
\newblock Mixing properties of {H}arris chains and autoregressive processes.
\newblock \emph{J. Appl. Probab.}, 23:\penalty0 880--892, 1986.

\bibitem[Baltagi(2001)]{b:bal-01}
B.~Baltagi.
\newblock \emph{Econometric Analysis of Panel Data}.
\newblock John Wiley and Sons, Chicester, 2nd edition, 2001.

\bibitem[Basrak et~al.(2002)Basrak, Davis, and Mikosch]{p:bas-03}
B.~Basrak, R.~Davis, and T.~Mikosch.
\newblock Regular variation of {GARCH} processes.
\newblock \emph{Stochastic Processes and their Applications}, 99:\penalty0
  95--115, 2002.

\bibitem[Bosq(1998)]{b:bos-98}
D.~Bosq.
\newblock \emph{Nonparametric Statistics for Stochastic Processes}.
\newblock Springer, New York, 1998.

\bibitem[Bousamma(1998)]{p:bou-98}
F.~Bousamma.
\newblock \emph{Ergodicit\'e, m\'elange et estimation dans les mod\`eles
  {GARCH}}.
\newblock PhD thesis, Paris 7, 1998.

\bibitem[Bradley(1996)]{p:bra-96}
R.~C. Bradley.
\newblock A covariance inequality under a two-part dependence assumption.
\newblock \emph{Statistics and Probability Letters}, 30:\penalty0 287--293,
  1996.

\bibitem[Chanda(1974)]{p:cha-74}
K.~C. Chanda.
\newblock Strong mixing properties of linear stochastic processes.
\newblock \emph{J. Appl. Prob.}, 11:\penalty0 401--408, 1974.

\bibitem[Cheng and Robinson(1991)]{p:che-rob-91}
B.~Cheng and P.~M. Robinson.
\newblock Density estimation in strongly dependent non-linear time series.
\newblock \emph{Statistica Sinica}, 1:\penalty0 335--359, 1991.

\bibitem[Cheng and Robinson(1994)]{p:che-rob-94}
B.~Cheng and P.~M. Robinson.
\newblock Semiparametric estimation from time series with long-range
  dependence.
\newblock \emph{Journal of Econometrics}, 94:\penalty0 335--353, 1994.

\bibitem[Cline and Pu(1999)]{p:cli-pu-99}
D.~Cline and H.~Pu.
\newblock Geometric ergodicity of nonlinear time series.
\newblock \emph{Statistica Sinica}, 9:\penalty0 1103--118, 1999.

\bibitem[Cs\"org\"o and Mielniczuk(1995)]{p:cso-95}
S.~Cs\"org\"o and J.~Mielniczuk.
\newblock Nonparametric regression under long-range dependent normal errors.
\newblock \emph{Ann. Statist.}, 23:\penalty0 1000--1014, 1995.

\bibitem[Cs\"org\"o and Mielniczuk(2001)]{p:cso-mie-01}
S.~Cs\"org\"o and J.~Mielniczuk.
\newblock The smoothing dichotomoy in random-design regression with long-memory
  based on moving averages.
\newblock \emph{Statistica Sinica}, 10:\penalty0 771--787, 2001.

\bibitem[Cs\"org\"o and Mielniczuk(1999)]{p:cso-mie-99}
S.~Cs\"org\"o and J.~Mielniczuk.
\newblock Random-design regression under long range dependence errors.
\newblock \emph{Bernoulli}, 5:\penalty0 209--224, 1999.

\bibitem[Dahlhaus and Feiler(2005)]{p:dah-fei-05}
R.~Dahlhaus and S.~Feiler.
\newblock Panel time series.
\newblock \emph{Preprint}, 2005.

\bibitem[Davidson(1994)]{b:dav-94}
J.~Davidson.
\newblock \emph{Stochastic Limit Theory}.
\newblock Oxford University Press, Oxford, 1994.

\bibitem[Doukhan(1994)]{b:dou-94}
P.~Doukhan.
\newblock \emph{Mixing, Properties and Examples}.
\newblock Springer, New York, 1994.

\bibitem[Estevas and Vieu(2003)]{p:vie-03}
G.~Estevas and P.~Vieu.
\newblock Nonparametric estimation under long memory dependence.
\newblock \emph{Nonparametric Statistics}, 15:\penalty0 535--551, 2003.

\bibitem[Geweke and Porter-Hudak(1983)]{p:gew-por-83}
J.~Geweke and S.~Porter-Hudak.
\newblock The estimation and application of long memory time series models.
\newblock \emph{Journal of Time Series Analysis}, 4:\penalty0 221--238, 1983.

\bibitem[Giraitis et~al.(1996)Giraitis, Koul, and Surgailis]{p:gir-kou-96}
L.~Giraitis, H.~L. Koul, and D.~Surgailis.
\newblock Asymptotic normality of regression estimators with long memory
  errors.
\newblock \emph{Statistics and Probability Letters}, 29:\penalty0 317--335,
  1996.

\bibitem[Giraitis et~al.(2000)Giraitis, Kokoskza, and Leipus]{p:gir-00}
L.~Giraitis, P.~Kokoskza, and R.~Leipus.
\newblock Stationary {ARCH} models: {D}ependence structure and central limit
  theorem.
\newblock \emph{Econometric Theory}, 16:\penalty0 3--22, 2000.

\bibitem[Gorodetskii(1977)]{p:gor-77}
V.~Gorodetskii.
\newblock On the strong mixing propery for linear sequences.
\newblock \emph{Theory of Probability and its Applications}, 22:\penalty0
  411--413, 1977.

\bibitem[Hall and Hart(1990{\natexlab{a}})]{p:hal-har-90}
P.~Hall and J.~Hart.
\newblock Convergence rates in density estimation for data from infinite-order
  moving average processes.
\newblock \emph{Probability Theory and Related Fields}, 87:\penalty0 253--274,
  1990{\natexlab{a}}.

\bibitem[Hall and Hart(1990{\natexlab{b}})]{p:hal-har-90b}
P.~Hall and J.~Hart.
\newblock Nonparametric regression with long-range dependence.
\newblock \emph{Stochastic Processes and their Applications}, 87:\penalty0
  339--351, 1990{\natexlab{b}}.

\bibitem[Hall and Heyde(1980)]{b:hal-hey-80}
P.~Hall and C.~Heyde.
\newblock \emph{{M}artingale {L}imit {T}heory and its {A}pplication}.
\newblock Academic Press, New York, 1980.

\bibitem[Hjellvik and Tj{\o}stheim(1999)]{p:tjo-99}
V.~Hjellvik and D.~Tj{\o}stheim.
\newblock Modelling panels of intercorrelated autoregressive time series
  models.
\newblock \emph{Biometrika}, 86:\penalty0 573--590, 1999.

\bibitem[Hjellvik et~al.(2004)Hjellvik, Chen, and Tj{\o}stheim]{p:tjo-04}
V.~Hjellvik, R.~Chen, and D.~Tj{\o}stheim.
\newblock Nonparametric estimation and testing in panels of intercorrelated
  time series models.
\newblock \emph{Journal of Time Series Analysis}, 25:\penalty0 831--872, 2004.

\bibitem[Johannes et~al.(2007)Johannes, Jun, and Subba~Rao]{p:joh-jun-sub-07}
J.~Johannes, M.~Jun, and S.~Subba~Rao.
\newblock Nonparametric estimation of spatio-temporal covariance functions.
\newblock \emph{Preprint}, 2007.

\bibitem[K{\"u}nsch(1987)]{p:kue-87}
H.~R. K{\"u}nsch.
\newblock Statistical aspects of self-similar processes.
\newblock \emph{Proceedings of the World Congress of the Bernoulli Society},
  1:\penalty0 67--74, 1987.

\bibitem[Linton and Mammen(2004)]{p:mam-04}
O.~B. Linton and E.~Mammen.
\newblock Estimating semiparametric {ARCH}($\infty$) models by kernel smoothing
  methods.
\newblock \emph{Econometrica}, 73:\penalty0 771--836, 2004.

\bibitem[Mammen et~al.(2005)Mammen, St{\o}ve, and Tj{\o}stheim]{p:mam-05}
E.~Mammen, B.~St{\o}ve, and D.~Tj{\o}stheim.
\newblock Nonparametric additive models for panels of time series.
\newblock \emph{Preprint}, 2005.

\bibitem[Masry and Tj{\o}stheim(1995)]{p:tjo-95}
E.~Masry and D.~Tj{\o}stheim.
\newblock Nonparametric estimation and identification of nonlinear {ARCH} time
  series.
\newblock \emph{Econometric Theory}, 11:\penalty0 258--289, 1995.

\bibitem[Mielniczuk(1997)]{p:mie-97}
J.~Mielniczuk.
\newblock On the asymptotic mean integrated squares error of kernel density
  estimator for dependent data.
\newblock \emph{Statistics and Probability Letters}, 34:\penalty0 53--58, 1997.

\bibitem[Rio(1993)]{p:rio-93}
E.~Rio.
\newblock Covariance inequalities for strongly mixing processes.
\newblock \emph{Ann. Inst. H. Poincar\'e Prob. Statist.}, 29:\penalty0
  587--597, 1993.

\bibitem[Robinson(1991)]{p:rob-91}
P.~M. Robinson.
\newblock Testing for strong serial correlation and dynamic conditional
  heteroskedasity in multiple regression.
\newblock \emph{Journal of Econometrics}, 47:\penalty0 67--78, 1991.

\bibitem[Robinson(1995)]{p:rob-95}
P.~M. Robinson.
\newblock Log-periodogram regression of time series with long range dependence.
\newblock \emph{Ann. Statist.}, 23:\penalty0 1048--1072, 1995.

\bibitem[Rosenblatt(1970)]{p:ros-70}
M.~Rosenblatt.
\newblock Density estimates and markov sequences.
\newblock In M.~Puri, editor, \emph{Non-parametric techniques in statistical
  inference}, pages 199--210. Cambridge University Press, London, 1970.

\bibitem[Scott(1992)]{b:S92}
D.~W. Scott.
\newblock \emph{Multivariate Density Estimation}.
\newblock Wiley, New York, 1992.

\bibitem[Subba~Rao(2006)]{p:sub-07}
S.~Subba~Rao.
\newblock A note on uniform convergence of an {ARCH}$(\infty)$ estimator.
\newblock \emph{Sankhya}, 68:\penalty0 600--620, 2006.

\bibitem[Subba~Rao(2007)]{p:sub-07b}
S.~Subba~Rao.
\newblock On mixing properties of {ARCH} and time-varying {ARCH} processes.
\newblock \emph{Preprint}, 2007.

\end{thebibliography}
